\documentclass[12pt]{amsart}
\usepackage{amssymb}
\usepackage{verbatim}
\usepackage{hyperref}
\usepackage{color}

\begin{document}

% define theorem environments
\newtheorem{theorem}{Theorem}    %[section]
\newtheorem{proposition}[theorem]{Proposition}
\newtheorem{conjecture}[theorem]{Conjecture}
\def\theconjecture{\unskip}
\newtheorem{corollary}[theorem]{Corollary}
\newtheorem{lemma}[theorem]{Lemma}
\newtheorem{sublemma}[theorem]{Sublemma}
\newtheorem{observation}[theorem]{Observation}
\newtheorem{remark}[theorem]{Remark}
\newtheorem{definition}[theorem]{Definition}
\theoremstyle{definition}
\newtheorem{notation}[theorem]{Notation}
\newtheorem{question}[theorem]{Question}
\newtheorem{questions}[theorem]{Questions}
\newtheorem{example}[theorem]{Example}
\newtheorem{problem}[theorem]{Problem}
\newtheorem{exercise}[theorem]{Exercise}

\numberwithin{theorem}{section} \numberwithin{theorem}{section}
\numberwithin{equation}{section}

\def\earrow{{\mathbf e}}
\def\rarrow{{\mathbf r}}
\def\uarrow{{\mathbf u}}
\def\varrow{{\mathbf V}}
\def\tpar{T_{\rm par}}
\def\apar{A_{\rm par}}

\def\reals{{\mathbb R}}
\def\torus{{\mathbb T}}
\def\heis{{\mathbb H}}
\def\integers{{\mathbb Z}}
\def\naturals{{\mathbb N}}
\def\complex{{\mathbb C}\/}
\def\distance{\operatorname{distance}\,}
\def\support{\operatorname{support}\,}
\def\dist{\operatorname{dist}\,}
\def\Span{\operatorname{span}\,}
\def\degree{\operatorname{degree}\,}
\def\kernel{\operatorname{kernel}\,}
\def\dim{\operatorname{dim}\,}
\def\codim{\operatorname{codim}}
\def\trace{\operatorname{trace\,}}
\def\Span{\operatorname{span}\,}
\def\dimension{\operatorname{dimension}\,}
\def\codimension{\operatorname{codimension}\,}
\def\nullspace{\scriptk}
\def\kernel{\operatorname{Ker}}
\def\ZZ{ {\mathbb Z} }
\def\p{\partial}
\def\rp{{ ^{-1} }}
\def\Re{\operatorname{Re\,} }
\def\Im{\operatorname{Im\,} }
\def\ov{\overline}
\def\eps{\varepsilon}
\def\lt{L^2}
\def\diver{\operatorname{div}}
\def\curl{\operatorname{curl}}
\def\etta{\eta}
\newcommand{\norm}[1]{ \|  #1 \|}
\def\expect{\mathbb E}
\def\bull{$\bullet$\ }
\def\C{\mathbb{C}}
\def\R{\mathbb{R}}
\def\Rn{{\mathbb{R}^n}}
\def\Sn{{{S}^{n-1}}}
\def\M{\mathbb{M}}
\def\N{\mathbb{N}}
\def\Q{{\mathbb{Q}}}
\def\Z{\mathbb{Z}}
\def\F{\mathcal{F}}
\def\L{\mathcal{L}}
\def\S{\mathcal{S}}
\def\supp{\operatorname{supp}}
\def\dist{\operatorname{dist}}
\def\essi{\operatornamewithlimits{ess\,inf}}
\def\esss{\operatornamewithlimits{ess\,sup}}
\def\xone{x_1}
\def\xtwo{x_2}
\def\xq{x_2+x_1^2}
\newcommand{\abr}[1]{ \langle  #1 \rangle}

\newcommand{\Norm}[1]{ \left\|  #1 \right\| }
\newcommand{\set}[1]{ \left\{ #1 \right\} }
\def\one{\mathbf 1}
\def\whole{\mathbf V}
\newcommand{\modulo}[2]{[#1]_{#2}}

\def\scriptf{{\mathcal F}}
\def\scriptg{{\mathcal G}}
\def\scriptm{{\mathcal M}}
\def\scriptb{{\mathcal B}}
\def\scriptc{{\mathcal C}}
\def\scriptt{{\mathcal T}}
\def\scripti{{\mathcal I}}
\def\scripte{{\mathcal E}}
\def\scriptv{{\mathcal V}}
\def\scriptw{{\mathcal W}}
\def\scriptu{{\mathcal U}}
\def\scriptS{{\mathcal S}}
\def\scripta{{\mathcal A}}
\def\scriptr{{\mathcal R}}
\def\scripto{{\mathcal O}}
\def\scripth{{\mathcal H}}
\def\scriptd{{\mathcal D}}
\def\scriptl{{\mathcal L}}
\def\scriptn{{\mathcal N}}
\def\scriptp{{\mathcal P}}
\def\scriptk{{\mathcal K}}
\def\frakv{{\mathfrak V}}
\allowdisplaybreaks
\begin{comment}
\def\scriptx{{\mathcal X}}
\def\scriptj{{\mathcal J}}
\def\scriptr{{\mathcal R}}
\def\scriptS{{\mathcal S}}
\def\scripta{{\mathcal A}}
\def\scriptk{{\mathcal K}}
\def\scriptp{{\mathcal P}}
\def\frakg{{\mathfrak g}}
\def\frakG{{\mathfrak G}}
\def\boldn{\mathbf N}
\end{comment}
%
%If a theorem-like environment should not be numbered,
%add * after \newtheorem, and delete the counter option such as [theorem].
\newtheorem*{remark0}{\indent\sc Remark}
%
%%%%% Proof %%%%%
\renewcommand{\proofname}{\indent\sc Proof.} % Name is small caps
%The following commands are available in the proof environment:
%\begin{proof}
%\end{proof}
%The end of a proof is marked with a square.
%%%%%%%%%%%%%%%%%%%%%%%%%%%%%%%%%%%%%%%%%

\title[Characterizations of BMO via endpoint bounds for commutators]
{Improvement to characterizations of BMO in terms of endpoint bounds for commutators and some counterexamples}
%title of paper and the running head option
\author{Yongming Wen and Huoxiong Wu}
%%%%%%%%%%%%%%% footnote %%%%%%%%%%%%%%%%
\subjclass[2010]{%2000 MSC numbers
42B20; 42B25; 42B35.
}
%In case \subjclass[2000] command is not effective
%(or the version of amsart.cls is old), write as follows instead:
%\renewcommand{\thefootnote}{\fnsymbol{footnote}}
%\footnote[0]{2000\textit{ Mathematics Subject Classification}.

%Primary 00; Secondary 00.}
%
\keywords{singular integrals, fractional integrals, Littlewood-Paley operators, variation operators, commutator, weak type endpoint estimate, BMO}
\thanks{Yongming Wen is supported by National Natural Science Foundation of China (Grant no. 12301119), Fujian Provincial Natural Science Foundation of China (Nos. 2025J01361), President's fund of Minnan Normal University (No. KJ2020020), Institute of Meteorological Big Data-Digital Fujian, Fujian Key Laboratory of Data Science and Statistics and Fujian Key Laboratory of Granular Computing and Applications (Minnan Normal University), China. Huoxiong Wu is supported by National Natural Science Foundation of China (Grant nos. 12271041, 12671114)}
%%%%%%%%%%%% Authors addresses %%%%%%%%%%%%%
\address{School of Mathematics and Statistics, Minnan Normal University, Zhangzhou 363000,  China} \email{wenyongmingxmu@163.com}
\address{School of Mathematical Sciences, Xiamen University, Xiamen 361005, China} \email{huoxwu@xmu.edu.cn}

%\address{}
%\email{}

%%%%%%%%%%%%%%%%%%%%%%%%%%%%%%%%%%%%%%%%%

\begin{abstract}
In this paper, we extend the technique of Lerner, Ombrosi and Rivera-R\'{i}os (Bull. Lond. Math. Soc. 51(1): 107--119, 2019) for the necessity of the strong-type boundedness of commutators to the endpoint case, and give characterizations of BMO via the endpoint bounds for commutators of singular integrals with rough kernels and related operators, which improve the conditions of Accomazzo (Israel J. Math. 228(2): 787--800, 2018), Guo et al. (J. Geom. Anal. 30(4): 3995-4035, 2020). We also give a characterization of BMO via the endpoint bounds for commutators of non-degenerate Calder\'{o}n-Zygmund operators by removing the real-valued restriction on pointwise multiplier, which extend the results of Hyt\"{o}nen (J. Math. Pures Appl. (9) 156: 351--391, 2021) to the endpoint case. Finally, we find counterexamples to illustrate that weighted endpoint estimate of commutator of fractional integral operator does not necessarily imply $b\in{\rm BMO}$, and two-weighted endpoint estimate of Bloom type for commutator of Hilbert transform does not hold, which are missing for a long time.
\end{abstract}

\maketitle

\section{Introduction and main results}
The aim of this paper is to characterize the class of symbols for which the commutators are bounded in the endpoint case. Let $b$ be a symbol, which is a locally integrable function on $\mathbb{R}^n$, and $T$ be a linear operator. We define the commutator $[b,T]$ by
\begin{align*}
[b,T]f:=T(bf)-bT(f)
\end{align*}
for suitable functions $f$. Let $R_j$ be the Riesz transform with $j=1,2,\cdots,n$ and $1<p<\infty$, the famous work of Coifman, Rochberg and Weiss \cite{CoRoWe} gave the following characterization:
\begin{align}\label{characterization}
\sum_{j=1}^n\|[b,R_j]\|_{L^p(\mathbb{R}^n)\rightarrow L^p(\mathbb{R}^n)}<\infty\Leftrightarrow b\in {\rm BMO}(\mathbb{R}^n),
\end{align}
where the function space BMO$(\mathbb{R}^n)$ consists of all locally integrable functions $b$ such that
\begin{align*}
\|b\|_{{\rm BMO}(\mathbb{R}^n)}:=\sup_{Q}\frac{1}{|Q|}\int_Q|b(x)-b_Q|dx<\infty,~ b_Q:=\frac{1}{|Q|}\int_Qb(x)dx.
\end{align*}
Consider homogeneous Calder\'{o}n-Zygmund kernel
\begin{align}\label{homogeneous CZ kernel}
K(x,y):=|x-y|^{-n}\Omega\Big(\frac{x-y}{|x-y|}\Big),
\end{align}
Jason \cite{Ja} and Uchiyama \cite{Uc} independently generalized \eqref{characterization} to commutators of singular integral operators $[b,T_\Omega]$ with $K$ satisfies certain smoothness condition, where
\begin{align}\label{convolution singular integral}
T_\Omega f(x):=p.v.\int_{\mathbb{R}^n}f(y)K(x,y)dy.
\end{align}
Since we have a wide array of tools at our disposal for addressing the $L^p$-boundedness of commutators, the primary difficulty in characterizing BMO spaces via the boundedness of commutators lies in the necessity. Over an extended timeframe, researchers have been unable to extend the BMO characterization to commutators of singular integral operators with rough kernels. Until recent years, Guo, Lian and Wu \cite{GLW} improved the condition of $\Omega$ by assuming there exists a nonempty open subset $E\subset \mathbb{S}^{n-1}$ such that
\begin{align}\label{upper and lower condition}
c\leq \Omega(x/|x|)\leq C
\end{align}
for all $x/|x|\in E$, where $0<c<C$ or $c<C<0$. Almost at the same time, by making use of the mean value method, Lerner, Ombrosi and Rivera-R\'{i}os \cite{LeOmRi} identified sufficient local positivity (lack of sign change in a nonempty open set) as a workable replacement of \eqref{upper and lower condition}. However, the method in \cite{LeOmRi} requires the pointwise multiplier to be real-valued. Taking advantage of approximate weak factorization, Hyt\"{o}nen \cite{Hy} perfectly solved the problems of relaxing the restriction on the pointwise multiplier and the sign-invariant condition on some open subset from $\mathbb{S}^{n-1}$ of $\Omega$. Recently, Li et al. \cite{LiLi} proposed the complex mean method to extend the results of Hyt\"{o}nen \cite{Hy} to higher order commutators.

For the endpoint case $p=1$, P\'{e}rez \cite{P} first proved that the commutators of Calder\'{o}n-Zygmund operators $[b,T]$ satisfy the following modular inequality when $b\in {\rm BMO}(\mathbb{R}^{n})$:
\begin{align}\label{LlogL}
\omega(\{x\in\mathbb{R}^{n}:|[b,T]f(x)|>t\})\lesssim\int_{\mathbb{R}^{n}}\Phi\Big(\frac{\|b\|_{{\rm BMO}(\mathbb{R}^{n})}|f(x)|}{t}\Big)\omega(x)dx,
~\omega\in A_1,
\end{align}
where $\Phi(t)=t\log(e+t)$. \eqref{LlogL} reflects that commutators are more singular than operators itself. When $\Omega\in Lip_\alpha(\mathbb{S}^{n-1})$ with $\alpha\in(0,1]$ and $b\in {\rm BMO}(\mathbb{R}^{n})$, \eqref{LlogL} also holds with $[b,T]$ replaced by $[b,T_\Omega]$.  Recently, Lan, Tao and Hu \cite{LaHuTa}  established the weighted weak type endpoint estimates for $[b,T_\Omega]$ with $\Omega\in L^q(\mathbb{S}^{n-1})(1<q\leq\infty)$. Let $\Omega$ be a function on $\mathbb{S}^{n-1}$. We say that $\Omega\in L(\log L)^2(\mathbb{S}^{n-1})$ if
\begin{align*}
\|\Omega\|_{L(\log L)^2(\mathbb{S}^{n-1})}:=\int_{\mathbb{S}^{n-1}}|\Omega(x')|\log^2(e+|\Omega(x')|)dx'<\infty.
\end{align*}
Hu and Tao \cite{HuTa} improved the aforementioned result to $\Omega\in L(\log L)^2(\mathbb{S}^{n-1})$. While for the necessity of the BMO$(\mathbb{R}^n)$ assumption for the endpoint estimates of commutators of singular integrals, the full picture is still very much incomplete. In \cite{A}, Accomazzo gave the following necessity criterion for $[b,T_\Omega]$.
\begin{theorem}{\rm(\cite{A})}\label{necessity criterion}
Let $b$ be a locally integrable function on $\mathbb{R}^n$ and $\Phi(t)=t\log(e+t)$. Let $\Omega\in Lip(\mathbb{S}^{n-1})$ be homogeneous of degree zero and $\int_{\mathbb{S}^{n-1}}\Omega=0$. Let $T_\Omega$ be defined as \eqref{convolution singular integral} and $T_\Omega^{\ast}$ be the adjoint operator. If there exists a constant $B$ such that for every measurable set $E$ and $t>0$, we have that
\begin{align*}
|\{x\in\mathbb{R}^n:|[b,T_\Omega](\chi_E)(x)|>t\}|\lesssim\int_{\mathbb{R}^n}\Phi\Big(\frac{B \chi_E(x)}{t}\Big)dx
\end{align*}
and
\begin{align*}
|\{x\in\mathbb{R}^n:|[b,T_\Omega^{\ast}](\chi_E)(x)|>t\}|\lesssim\int_{\mathbb{R}^n}\Phi\Big(\frac{B \chi_E(x)}{t}\Big)dx.
\end{align*}
Then $b\in {\rm BMO(\mathbb{R}^n)}$.
\end{theorem}
Obviously, if $T$ is the Riesz transform, the condition on the adjoint would be free. Almost at the same time, Guo, Lian and Wu \cite{GLW} also proved the necessity for weighted endpoint estimates of $[b,T_\Omega]$ without the adjoint condition, however, the technique still relies upon the regularity of the kernel. Hence, it is very natural to ask the following questions.

\textbf{Question 1:} Inspired by the results in \cite{Hy,LeOmRi}, is it possible to remove the smooth assumption of $\Omega$ in \cite{A,GLW}?

Our first result is concerned with the characterization of ${\rm BMO}$ in terms of endpoint bound for commutator of singular integral with rough kernel.

\begin{theorem}\label{characterization endpoint estimate of commutator of CZO with b real}
Let $T_\Omega$ be defined by \eqref{convolution singular integral}, $\Phi(t)=t\log(e+t)$, $b\in L_{loc}^1(\mathbb{R}^n)$ be a real-valued function. Assume that $\Omega$ does not change
sign and is not equivalent to zero on some open subset of $S^{n-1}$, $\Omega(\lambda x)=\Omega(x)$ for any $\lambda>0$ and $x\in\mathbb{R}^n$, $\int_{\mathbb{S}^{n-1}}\Omega=0$ and $\Omega\in L(\log L)^2(\mathbb{S}^{n-1})$. Then the following two statements are equivalent.
\begin{itemize}
\item [(1)] $b\in{\rm BMO}(\mathbb{R}^n)$;
\item [(2)] there exists constants $c,\alpha>0$ such that for any $t>0$,
\begin{align*}
|\{x\in\mathbb{R}^n:|[b,T_\Omega]f(x)|>t\}|\leq c\int_{\mathbb{R}^n}\Phi\Big(\frac{\alpha|f(x)|}{t}\Big)dx.
\end{align*}
\end{itemize}
\end{theorem}

We also establish the characterizations of ${\rm BMO}$ in terms of endpoint bounds for $\mathcal{V}_{\rho}(\mathcal{T}_{\Omega,b}f)$, $\mathcal{V}_\rho((\Phi\star f)_b)$, $u_{\Omega,b}$, $u_{S,b}^\varrho$ and $g_{\lambda,b}^{\ast,\varrho}$. Before we state these characterizations, we first recall the definitions of these commutators.
\begin{definition}
$(1)$ Let $K$ be defined as \eqref{homogeneous CZ kernel}, where $\Omega$ is a homogeneous function of degree zero on $\mathbb{R}^n$, $\Omega\in L^1(\mathbb{S}^{n-1})$ and $\int_{\mathbb{S}^{n-1}}\Omega=0$. Let $b\in L^1_{\rm loc}(\mathbb{R}^n)$, $\mathcal{T}_{\Omega,b}:=\{T_{\varepsilon,\Omega,b}\}_{\varepsilon>0}$, where
\begin{equation*}
T_{\varepsilon,\Omega,b}(f)(x):=\int_{|x-y|>\varepsilon}(b(x)-b(y))K(x,y)f(y)dy.
\end{equation*}
Define the $\rho$-variation operator for commutator of singular integral as follows
$$\mathcal{V}_\rho(\mathcal{T}_{\Omega,b}f)(x):=\sup_{\varepsilon_i\downarrow0}\Big(\sum_{i=1}^{\infty}|
T_{\varepsilon_{i+1},\Omega,b}f(x)-T_{\varepsilon_i,\Omega,b}f(x)|^\rho\Big)^{{1}/{\rho}},~\rho>2,$$
where the supremum is taken over all sequence $\{\varepsilon_i\}$ decreasing to zero.

$(2)$ Let $b\in L^1_{\rm loc}(\mathbb{R}^n)$, $\phi\in\mathcal{S}(\mathbb{R}^n)$ satisfy $\int_{\mathbb{R}^n}\phi(x)dx=1$, where $\mathcal{S}(\mathbb{R}^n)$ is the space of Schwartz functions. Consider the family of commutators associated with approximate identities, which is defined by
\begin{equation*}
(\Phi\star f)_b(x):=\{b(x)(\phi_t\ast f)(x)-\phi_t\ast (bf)(x)\}_{t>0},
\end{equation*}
where $\phi_t(x)=t^{-n}\phi(x/t)$,
$$b(x)(\phi_t\ast f)(x)-\phi_t\ast (bf)(x):=\int_{\mathbb{R}^n}\frac{1}{t^n}\phi\Big(\frac{x-y}{t}\Big)(b(x)-b(y))f(y)dy.$$
The $\rho$-variation operator associated with commutators of approximate identities is defined by
\begin{align*}
\mathcal{V}_\rho((\Phi\star f)_b)(x)&:=\sup_{\{t_k\}\downarrow0}\Big(\sum_{k\geq1}
\Big|b(x)(\phi_{t_k}\ast f)(x)-\phi_{t_k}\ast (bf)(x)\\
&\qquad\qquad-b(x)(\phi_{t_{k+1}}\ast f)(x)+\phi_{t_{k+1}}\ast (bf)(x)
\Big|^\rho\Big)^{1/\rho},~\rho>2.
\end{align*}
\end{definition}

\begin{definition}
Let $\Omega$ be a homogeneous function of degree zero on $\mathbb{R}^n$, $\Omega\in L^1(\mathbb{S}^{n-1})$ and $\int_{\mathbb{S}^{n-1}}\Omega=0$. Let $b\in L^1_{\rm loc}(\mathbb{R}^n)$, $\varrho>n/2$ and $\lambda>2$. The commutator of Marcinkiewicz type integral $u_{\Omega,b}$, the commutator of parameterized area integral $u_{S,b}^\varrho$ and the commutator of parameterized Littlewood-Paley $g_{\lambda,b}^{\ast,\varrho}$ are defined by
\begin{align*}
u_{\Omega,b}f(x):=\Big(\int_0^\infty\Big|\int_{|x-y|\leq t}(b(x)-b(y))\frac{\Omega(x-y)}{|x-y|^{n-1}}f(y)dy\Big|^2\frac{dt}{t^3}\Big)^{1/2},
\end{align*}
\begin{align*}
u_{S,b}^\varrho(f)(x):=\Big(\iint_{\Gamma(x)}\Big|t^{-\varrho}\int_{|y-z|<t}
\frac{\Omega(y-z)}{|y-z|^{n-\varrho}}(b(x)-b(z))
f(z)dz\Big|^2\frac{dydt}{t^{n+1}}\Big)^{1/2},
\end{align*}
where $\Gamma(x):=\{(y,t)\in\mathbb{R}_+^{n+1}:|x-y|<t\}$,
and
\begin{align*}
&g_{\lambda,b}^{\ast,\varrho}(f)(x)\\
&\quad:=\Big(\iint_{\mathbb{R}_+^{n+1}}\Big(\frac{t}{t+|x-y|}\Big)^{\lambda n}\Big|t^{-\varrho}\int_{|y-z|<t}\frac{\Omega(y-z)}{|y-z|^{n-\varrho}}(b(x)-b(z))
f(z)dz\Big|^2\frac{dydt}{t^{n+1}}\Big)^{1/2},
\end{align*}
respectively.
\end{definition}

Let $\{T_t\}_{t>0}$ be a family of operators that are bounded on certain function spaces, where $\lim_{t\rightarrow0}T_t f$ converges in some sense. Variation operator serves as a substitute for the square function to quantify the rate of convergence of $\{T_t\}_{t>0}$. If variation operator belongs to some function spaces such as $L^p$ or $L^{p,\infty}$, the corresponding sequence converges almost everywhere with no extra conditions required, which are not available in certain ergodic models, see \cite{Bo}. Variation operator was first investigated by L\'{e}pingle \cite{Le} to improve the classical Doob maximal inequality. Later on, Bourgain \cite{Bo} studied variational estimates for ergodic averages of a dynamic system. Since then, numerous researchers have studied variation of several families of operators across probability, ergodic theory and harmonic analysis, for instance, see \cite{AJS,CaJRW1,CaJRW2,Jo,JR} for prior work. Among which, Campbell et al. \cite{CaJRW1,CaJRW2} first established the $L^p$-bounds and weak type $(1,1)$ bounds for variation of truncated singular integrals. Liu, Wu \cite{LW} and Betancor et al. \cite{BFHR} independently established $L^p$-boundedness for variation operators of commutators of singular integrals with smooth kernels. Chen et al. \cite{ChDiHoLi} extended the results in \cite{LW} to variation operators of commutators of singular integrals with rough kernels. The weighted endpoint estimates for $\mathcal{V}_\rho(\mathcal{T}_{\Omega,b}f)$ are established recently in \cite{GuoWenWuYang,WenWuXue}, here, we focus on the reverse results, which are given below.
\begin{theorem}\label{characterization endpoint estimate of commutator of VCZO with b real}
Let $\rho>2$, $\Phi(t)=t\log(e+t)$, $b\in L_{loc}^1(\mathbb{R}^n)$ be a real-valued function and $w\in A_1$. Suppose that $\Omega$ is not identically zero, $\Omega(\lambda x)=\Omega(x)$ for any $\lambda>0$ and $x\in\mathbb{R}^n$, $\int_{\mathbb{S}^{n-1}}\Omega=0$ and $\Omega\in {\rm Lip}_\alpha(S^{n-1})$ with $0<\alpha\le 1$. Then the following two statements are equivalent.
\begin{itemize}
\item [(1)] $b\in{\rm BMO}(\mathbb{R}^n)$;
\item [(2)] there exists constants $c,\alpha>0$ such that for any $t>0$,
\begin{align*}
w(\{x\in\mathbb{R}^n:\mathcal{V}_{\rho}(\mathcal{T}_{\Omega,b}f)(x)>t\})\leq c\int_{\mathbb{R}^n}\Phi\Big(\frac{\alpha|f(x)|}{t}\Big)w(x)dx.
\end{align*}
\end{itemize}
\end{theorem}

\begin{theorem}\label{characterization endpoint estimate of commutator of V approximate}
Let $\rho>2$, $\Phi(t)=t\log(e+t)$, $\phi\in \mathcal S(\mathbb R^n)$ with $\int_{\mathbb{R}^n}\phi(x)dx=1$, $w\in A_1$. If $\phi$ or $b$ is a real-valued function, then for $\rho>2$, the following statements are equivalent:
\begin{itemize}
\item [(1)] $b\in{\rm BMO}(\mathbb{R}^n)$;
\item [(2)] there exists constants $c,\alpha>0$ such that for any $t>0$,
\begin{align*}
w(\{x\in\mathbb{R}^n:\mathcal{V}_\rho((\Phi\star f)_b)(x)>t\})\leq c\int_{\mathbb{R}^n}\Phi\Big(\frac{\alpha|f(x)|}{t}\Big)w(x)dx.
\end{align*}
\end{itemize}
\end{theorem}

\begin{theorem}\label{characterization endpoint estimate of commutator of Marcinkiewicz}
Let $\rho>2$, $\Phi(t)=t\log(e+t)$, $\phi\in \mathcal S(\mathbb R^n)$ with $\int_{\mathbb{R}^n}\phi(x)dx=1$, $w\in A_1$. If $\phi$ or $b$ is a real-valued function, then for $\rho>2$, the following statements are equivalent:
\begin{itemize}
\item [(1)] $b\in{\rm BMO}(\mathbb{R}^n)$;
\item [(2)] there exists constants $c,\alpha>0$ such that for any $t>0$,
\begin{align*}
w(\{x\in\mathbb{R}^n:\mathcal{V}_\rho((\Phi\star f)_b)(x)>t\})\leq c\int_{\mathbb{R}^n}\Phi\Big(\frac{\alpha|f(x)|}{t}\Big)w(x)dx.
\end{align*}
  \end{itemize}
\end{theorem}

It is well-known that Littlewood-Paley operators play very important roles in harmonic analysis and PDE, see \cite{ChWiWo,Ke}. The higher dimension Marcinkiewicz integral operator was first introduced by Stein \cite{St}, who proved $u_\Omega$ is strong type $(p,p)$ and weak type $(1,1)$. The first result bout the commutators of the Marcinkiewicz integrals was given by Torchinsky and Wang \cite{TW}, in which the authors established weighted $L^p$ estimates for $u_{\Omega,b}$. Zhang \cite{Zhang} considered the weighted endpoint estimates for $u_{\Omega,b}$. As generalizations of area integral and Littlewood-Paley $g_\lambda^{\ast}$ function, parameterized Littlewood-Paley and parameterized area integral were first introduced by Sakamoto and Yabuta \cite{SaYa}, they gave $L^p$ boundedness for parameterized Littlewood-Paley and parameterized area integral, while the weak type $(1,1)$ estimates for parameterized Littlewood-Paley and parameterized area integral were later proved by Ding and Xue \cite{DiXu0}. Subsequently, Ding and Xue \cite{DiXu} established weighted strong type and weak type estimates for $u_{S,b}^\varrho$ and $g_{\lambda,b}^{\ast,\varrho}$.

Continuing the study in \cite{DiXu,Zhang}, in which the authors gave the weighted endpoint estimates for $u_{\Omega,b}$, $u_{S,b}^\varrho$ and $g_{\lambda,b}^{\ast,\varrho}$, we aim to show that $b\in{\rm BMO(\mathbb{R}^n)}$ is also necessary for the weighted endpoint estimates of $u_{\Omega,b}$, $u_{S,b}^\varrho$ and $g_{\lambda,b}^{\ast,\varrho}$.
\begin{theorem}\label{characterization endpoint estimate of commutator of Marcinkiewicz with b real}
Let $\Phi(t)=t\log(e+t)$, $b\in L_{loc}^1(\mathbb{R}^n)$ be a real-valued function and $w\in A_1$. Assume that $\Omega$ does not change sign and is not equivalent to zero on some open subset of $S^{n-1}$, $\Omega(\lambda x)=\Omega(x)$ for any $\lambda>0$ and $x\in\mathbb{R}^n$, $\int_{\mathbb{S}^{n-1}}\Omega=0$, $\Omega\in L^\infty(\mathbb{S}^{n-1})$ and
\begin{align*}
\int_0^1\frac{\omega_\infty(\delta)}{\delta}\log\frac{1}{\delta}d\delta<\infty,
\end{align*}
where
\begin{align*}
\omega_\infty(\delta):=\sup_{\|\rho\|<\delta}\sup_{x'\in\mathbb{S}^{n-1}}|\Omega(\rho x')-\Omega(x')|,
\end{align*}
and $\|\rho\|:=\sup_{x'\in\mathbb{S}^{n-1}}|\rho x'-x'|$. Then the following two statements are equivalent.
\begin{itemize}
\item [(1)] $b\in{\rm BMO}(\mathbb{R}^n)$;
\item [(2)] there exists constants $c,\alpha>0$ such that for any $t>0$,
\begin{align*}
w(\{x\in\mathbb{R}^n:u_{\Omega,b}(f)(x)>t\})\leq c\int_{\mathbb{R}^n}\Phi\Big(\frac{\alpha|f(x)|}{t}\Big)w(x)dx.
\end{align*}
\end{itemize}
\end{theorem}

\begin{theorem}\label{characterization endpoint estimate of commutator of parameter with b real}
Let $T$ be $u_{S,b}^\varrho$ or $u_{\lambda,b}^{\ast,\varrho}$. Let $\varrho>n/2$, $\lambda>2$, $\Phi(t)=t\log(e+t)$, $b\in L_{loc}^1(\mathbb{R}^n)$ be a real-valued function and $w\in A_1$. Assume that $\Omega$ does not change sign and is not equivalent to zero on some open subset of $\mathbb{S}^{n-1}$, $\Omega(\lambda x)=\Omega(x)$ for any $\lambda>0$ and $x\in\mathbb{R}^n$, $\int_{\mathbb{S}^{n-1}}\Omega=0$, $\Omega\in L^2(\mathbb{S}^{n-1})$ and
\begin{align*}
\int_0^1\frac{\omega_2(\delta)}{\delta}(1+|\log \delta|)^\sigma d\delta<\infty,~\sigma>2,
\end{align*}
where
\begin{align*}
\omega_2(\delta):=\sup_{\|\rho\|<\delta}\Big(\int_{\mathbb{S}^{n-1}}|\Omega(\rho x')-\Omega(x')|^2d\sigma(x')\Big)^{1/2},
\end{align*}
and $\|\rho\|:=\sup_{x'\in\mathbb{S}^{n-1}}|\rho x'-x'|$. Then the following two statements are equivalent.
\begin{itemize}
\item [(1)] $b\in{\rm BMO}(\mathbb{R}^n)$;
\item [(2)] there exists constants $c,\alpha>0$ such that for any $t>0$,
\begin{align*}
w(\{x\in\mathbb{R}^n:T(f)(x)>t\})\leq c\int_{\mathbb{R}^n}\Phi\Big(\frac{\alpha|f(x)|}{t}\Big)w(x)dx.
\end{align*}
\end{itemize}
\end{theorem}

\begin{remark}
From the proofs of Theorem \ref{characterization endpoint estimate of commutator of CZO with b real}, Theorems \ref{characterization endpoint estimate of commutator of VCZO with b real}-\ref{characterization endpoint estimate of commutator of parameter with b real}, we emphasize that the process of $(2)\Rightarrow(1)$ does not require any regularity of the kernel.
\end{remark}

Recall that we call that $K$ is an $\omega$-Calder\'{o}n-Zygmund kernel if $K$ satisfies
\begin{align*}
|K(x,y)|\leq\frac{c_K}{|x-y|^n},~x\neq y,
\end{align*}
\begin{align*}
|K(x,y)-K(x',y)|+|K(y,x)-K(y,x')|\leq\frac{1}{|x-y|^n}\omega\Big(\frac{|x-x'|}{|x-y|}\Big),
\end{align*}
whenever $|x-x'|<|x-y|/2$, where $\omega:[0,1)\rightarrow[0,\infty)$ is an increasing subadditive continuous function satisfying $\omega(0)=0$ and the Dini-condition
\begin{align*}
\int_0^1\omega(t)\frac{dt}{t}<\infty.
\end{align*}

\begin{definition}{\rm(\cite{Hy})}\label{non-degenerate CZ}
We say that $K$ is a non-degenerate Calder\'{o}n-Zygmund kernel, if (at least) one of the following two conditions holds:
\begin{itemize}
\item [(1)] $K$ is an $\omega$-Calder\'{o}n-Zygmund kernel and for every $y\in\mathbb{R}^n$ and $r>0$, there exists $x\in B(y,r)^c$ with $|K(x,y)|\geq1/(c_0r^n)$.
\item [(2)]  $K$ is a homogeneous Calder\'{o}n-Zygmund kernel with $\Omega\in L^1(\mathbb{S}^{n-1})$. In particular, there exists a Lebesgue point $\theta_0\in\mathbb{S}^{n-1}$ of $\Omega$ such that $\Omega(\theta_0)\neq0$.
\end{itemize}
\end{definition}

The following theorem extends the result of Hyt\"{o}nen \cite{Hy} to the endpoint case. Compared with the theorems we establish above, we remove the restriction of real-valued condition on $b$ for commutators of non-degenerate Calder\'{o}n-Zygmund operators.

\begin{theorem}\label{necesity endpoint estimate of commutator of non-degenerate CZO}
Let $T$ be an non-degenerate Calder\'{o}n-Zygmund operator. Let $\Phi(t)=t\log(e+t)$, $b\in L_{loc}^1(\mathbb{R}^n)$ and $w\in A_1$. If
\begin{align*}
w(\{x\in\mathbb{R}^n:|[b,T]f(x)|>t\})\leq c\int_{\mathbb{R}^n}\Phi\Big(\frac{\alpha|f(x)|}{t}\Big)w(x)dx,
\end{align*}
then $b\in{\rm BMO}(\mathbb{R}^n)$.
\end{theorem}
\begin{remark}
We believe that the technique used in Theorem \ref{necesity endpoint estimate of commutator of non-degenerate CZO} can also applied to Theorems \ref{characterization endpoint estimate of commutator of VCZO with b real}-\ref{characterization endpoint estimate of commutator of parameter with b real}, so the restriction of real-valued condition on $b$ can also be relaxed. However, compared with the proofs of Theorems \ref{characterization endpoint estimate of commutator of VCZO with b real}-\ref{characterization endpoint estimate of commutator of parameter with b real}, the proof of \ref{necesity endpoint estimate of commutator of non-degenerate CZO} are more complex. Therefore, it is tedious to prove similar results for $\mathcal{V}_{\rho}(\mathcal{T}_{\Omega,b}f)$, $\mathcal{V}_\rho((\Phi\star f)_b)$, $u_{\Omega,b}$, $u_{S,b}^\varrho$, $g_{\lambda,b}^{\ast,\varrho}$, we leave details for interested readers.
\end{remark}

By Theorem \ref{necesity endpoint estimate of commutator of non-degenerate CZO} and the endpoint estimates of non-degenerate Calder\'{o}n-Zygmund operators in \cite{HuTa,LeOmRi2017}, we have the following two corollaries.
\begin{corollary}\label{characterization endpoint estimate of commutator of non-degenerate CZO1}
Let $T$ be an non-degenerate $\omega$-Calder\'{o}n-Zygmund operator. Let $\Phi(t)=t\log(e+t)$, $b\in L_{loc}^1(\mathbb{R}^n)$ and $w\in A_1$. Then the following two statements are equivalent.
\begin{itemize}
\item [(1)] $b\in{\rm BMO}(\mathbb{R}^n)$.
\item [(2)] there exists constants $c,\alpha>0$ such that for any $t>0$,
\begin{align*}
w(\{x\in\mathbb{R}^n:|[b,T]f(x)|>t\})\leq c\int_{\mathbb{R}^n}\Phi\Big(\frac{\alpha|f(x)|}{t}\Big)w(x)dx.
\end{align*}
\end{itemize}
\end{corollary}

\begin{corollary}\label{characterization endpoint estimate of commutator of non-degenerate CZO2}
Let $T$ be an non-degenerate homogeneous Calder\'{o}n-Zygmund operator with $\Omega(\lambda x)=\Omega(x)$ for any $\lambda>0$, $\int_{\mathbb{S}^{n-1}}\Omega=0$ and $\Omega\in L(\log L)^2(\mathbb{S}^{n-1})$. Let $\Phi(t)=t\log(e+t)$, $b\in L_{loc}^1(\mathbb{R}^n)$. Then the following two statements are equivalent.
\begin{itemize}
\item [(1)] $b\in{\rm BMO}(\mathbb{R}^n)$.
\item [(2)] there exists constants $c,\alpha>0$ such that for any $t>0$,
\begin{align*}
|\{x\in\mathbb{R}^n:|[b,T]f(x)|>t\}|\leq c\int_{\mathbb{R}^n}\Phi\Big(\frac{\alpha|f(x)|}{t}\Big)dx.
\end{align*}
\end{itemize}
\end{corollary}

Let $[b,I_\alpha]$ be the commutator of fractional integral, which is defined by
\begin{align*}
[b,I_\alpha]f(x):=\int_{\mathbb{R}^n}\frac{(b(x)-b(y))f(y)}{|x-y|^{n-\alpha}}dy,~0<\alpha<n.
\end{align*}
Chanillo \cite{Ch} first established $[b,I_\alpha]:L^p(\mathbb{R}^n)\rightarrow L^q(\mathbb{R}^n)\Leftrightarrow b\in{\rm BMO}(\mathbb{R}^n)$. Segovia and Torrea \cite{ST} extended this result to the weighted case. Cruz-Uribe, Fiorenza \cite{CrFi2003,CrFi2007} and Ding, Lu, Zhang \cite{DiLuZh} independently established the following endpoint estimates of $[b,I_\alpha]$,
\begin{align*}
b\in{\rm BMO}(\mathbb{R}^n)\Rightarrow w(\{x\in\mathbb{R}^n:|[b,I_\alpha]f(x)|>t\})\lesssim \Psi\Big(\int_{\mathbb{R}^n}\Phi\Big(\frac{c|f(x)|}{t}\Big)\Theta(w(x))dx\Big).
\end{align*}
However, the reverse result is a long-standing problem that has not been solved.

\textbf{Question 2:} Can we establish
\begin{align*}
\|b\|_{{\rm BMO}(\mathbb{R}^n)}\sim w(\{x\in\mathbb{R}^n:|[b,I_\alpha]f(x)|>t\})\lesssim \Psi\Big(\int_{\mathbb{R}^n}\Phi\Big(\frac{|f(x)|}{t}\Big)\Theta(w(x))dx\Big)?
\end{align*}

Surprisingly, we give a negative answer. The weighted endpoint estimate of $[b,I_\alpha]$ does not necessarily imply $b\in{\rm BMO}(\mathbb{R}^n)$.
\begin{theorem}\label{counterexample}
Let $0<\alpha<n$, $q=n/(n-\alpha)$, $\Phi(t)=t\log(e+t)$, $\Psi(t)=[t\log(e+t^{\alpha/n})]^q$, $\Theta(t)=t^{1/q}\log(e+t^{-\alpha/n})$ and $w\in A_1$. There exists a function $b\not\in{\rm BMO}(\mathbb{R}^n)$ such that
\begin{align*}
w(\{x\in\mathbb{R}^n:|[b,I_\alpha]f(x)|>t\})\leq C_{w,n,\alpha}\Psi\Big(\int_{\mathbb{R}^n}\Phi\Big(\frac{|f(x)|}{t}\Big)\Theta(w(x))dx\Big).
\end{align*}
\end{theorem}

On the other hand, in one dimension, Bloom \cite{Bloom} extended \eqref{characterization} to the following weighted case:
\begin{align}\label{Bloom Hilbert}
\|[b,H]\|_{L^p(\mu)\rightarrow L^p(\lambda)}\sim\|b\|_{{\rm BMO}_\eta(\mathbb{R})},
\end{align}
where $\mu,\lambda\in A_p$, $\eta=(\mu/\lambda)^{1/p}$, $H$ is the Hilbert transform
\begin{align*}
Hf(x)=p.v.\int_{\mathbb{R}}\frac{f(y)}{x-y}dy,
\end{align*}
and the weighted ${\rm BMO}$ space ${\rm BMO}_\eta(\mathbb{R}^n)$ is defined by
\begin{align*}
{\rm BMO}_\eta(\mathbb{R}^n):=\Big\{b\in L_{loc}^1(\mathbb{R}^n):\|b\|_{{\rm BMO}_\eta(\mathbb{R}^n)}:=\sup_{B\subset\mathbb{R}^n}\frac{1}{\eta(B)}\int_B|b(x)-b_B|dx<\infty\Big\}.
\end{align*}
We refer to this two-weighted estimate as Bloom-type estimate. Until 2016, Holmes, Lacey and Wick \cite{HLW} extended \eqref{Bloom Hilbert} to Riesz transform. Lerner et al. \cite{LeOmRi} established the Bloom-type characterization of the two-weighted boundedness of commutators of singular integrals. Hyt\"{o}nen \cite{Hy} further extended the results in \cite{LeOmRi} to non-degenerate $w$-Calder\'{o}n-Zygmund operator. Combined with the results in \cite{Li} and \cite{Hy}, the following Bloom-type characterization of the two-weighted boundedness of commutators of non-degenerate homogeneous Calder\'{o}n-Zygmund operators is established
\begin{align*}
\|[b,T_\Omega]\|_{L^p(\mu)\rightarrow L^p(\lambda)}\sim\|b\|_{{\rm BMO}_\eta(\mathbb{R}^n)},~\Omega\in L^\infty(\mathbb{S}^{n-1}).
\end{align*}
For more results about this topic, see \cite{AMR,DiLu1999,DuGoKuLiWiYa,HuLiWa,LiMaVu,ST1993} and the references therein. Up to now, we have known nothing about two-weighted endpoint estimates of Bloom type for commutators of singular integrals.

\textbf{Question 3:} Is it possible to establish two-weighted endpoint estimates of Bloom type for commutators of singular integrals? And then give a characterization of weighted ${\rm BMO(\mathbb{R}^n)}$ space via two-weighted endpoint estimates of Bloom type for commutators of singular integrals?

The following theorem gives a negative answer.
\begin{theorem}\label{counterexample1}
Let $H$ be the Hilbert transform and $\Phi(t)=t\log(e+t)$. Then there exist weights $\lambda,\mu\in A_1$, functions $f\in L_c^\infty(\mathbb{R})$, $b\in{\rm BMO}_\eta(\mathbb{R})$ with $\eta=\mu/\lambda$ and $t>0$ such that
\begin{align*}
\lambda(\{x\in\mathbb{R}:|[b,H]f(x)|>t\})\geq C\int_{\mathbb{R}}\Phi\Big(\frac{|f(x)|}{t}\Big)\mu(x)dx.
\end{align*}
\end{theorem}

We organize the rest of the paper as follows. In Section 2, we will prove Theorem \ref{characterization endpoint estimate of commutator of CZO with b real}, Theorems \ref{characterization endpoint estimate of commutator of VCZO with b real}-\ref{characterization endpoint estimate of commutator of parameter with b real}. Theorem \ref{necesity endpoint estimate of commutator of non-degenerate CZO} will be proved in Section 3. The proof of Theorem \ref{counterexample} and Theorem \ref{counterexample1} will be given in Section 4 and Section 5, respectively.

Throughout the rest of the paper, we denote a positive constant by $C$, which may change at each occurrence. We also denote $f\lesssim g$, $f\thicksim g$ if $f\leq Cg$ and $f\lesssim g \lesssim f$, respectively. For any ball $B:=B(x_B,r_B)\subset \mathbb{R}^n$ and $\sigma>0$, $\chi_B$ represents the characteristic function of $B$, $\sigma B$ means $B(x_B,\sigma r_B)$, $w_B:=|B|^{-1}\int_Bw:=|B|^{-1}w(B)$.

\section{Characterizations of the weighted endpoint estimates of commutators of singular integral operators and related operators}
In this section, we prove Theorems \ref{characterization endpoint estimate of commutator of CZO with b real}, \ref{characterization endpoint estimate of commutator of VCZO with b real}-\ref{characterization endpoint estimate of commutator of parameter with b real}. The main tools in this section are median value and local mean oscillation of a real-valued measurable function, which were introduced in \cite{LeOmRi} to solve the
necessity of ${\rm BMO}$ for the strong type estimate of commutators. Now, we will extend this technique to the endpoint case. We recall some definitions and lemmas first.

\begin{definition}
A weight is a nonnegative and locally integrable function on $\mathbb{R}^n$.
We say that $w\in A_{1}$ if
$$Mw(x)\leq cw(x),$$
where the smallest constant $c$ is refer to as $[w]_{A_1}$.
\end{definition}

\begin{definition} {\rm (\cite{LeOmRi})}
By a median value of a real-valued measurable function $f$ over a measure set $E$ of positive finite measure, we mean a possibly non-unique, real number $m_f(E)$ such that
$$\max(|\{x\in E: f(x)>m_f(E)\}|,\,\,|\{x\in E: f(x)<m_f(E)\}|)\leq|E|/2.$$
\end{definition}
\begin{definition}{\rm(\cite{LeOmRi})}
Let $0<\lambda<1$. Given a measurable function $f$ on $\mathbb{R}^n$ and a cube $Q$, the local mean oscillation of $f$ on $Q$ is
defined by
\begin{align*}
w_\lambda(f;Q)=\inf_{c\in\mathbb{R}}((f-c)\chi_Q)^{\ast}(\lambda|Q|),
\end{align*}
where $f^{\ast}$ denotes the non-increasing rearrangement of $f$.
\end{definition}

It was proved in \cite{LeOmRi} that
\begin{align}\label{transform to local mean oscillation}
\|f\|_{{\rm BMO(\mathbb{R}^n)}}\lesssim w_\lambda(b;Q),~0<\lambda\leq2^{-n-2}.
\end{align}

\begin{lemma}{\rm(\cite{LeOmRi})}\label{mean value method}
There exist $0<\epsilon_0,\xi_0<1$ and $k_0>1$ depending only on $\Omega$ and $n$ such that the following holds. For every cube $Q\subset\mathbb{R}^n$, there exist measurable sets $E\subset Q$, $F\subset k_0Q$ with $|E|\sim|F|\sim|Q|$ and $G\subset E\times F$ with $|G|\geq\xi_0|Q|^2$ such that
\begin{itemize}
\item [(i)] $w_{2^{-n-2}}(b;Q)\leq|b(x)-b(y)|$ for all $(x,y)\in E\times F$;
\item [(ii)] $\Omega\big(\frac{x-y}{|x-y|}\big)$ and $b(x)-b(y)$ do not change sign in $E\times F$;
\item [(iii)] $|\Omega\big(\frac{x-y}{|x-y|}\big)|\geq\epsilon_0$ for all $(x,y)\in G$.
\end{itemize}
\end{lemma}

\begin{proof}[Proof of Theorem \ref{characterization endpoint estimate of commutator of CZO with b real}]
$(1)\Rightarrow(2)$ is proved in \cite{LeOmRi2017}, we only need to prove $(2)\Rightarrow(1)$. Given a cube $Q$, let $E$, $F$, $G$ be the sets stated in Lemma \ref{mean value method}. Let $f=\chi_F$. By Lemma \ref{mean value method}, we have
\begin{align*}
|[b,T_{\Omega}](f)(x)\chi_E(x)|&=\Big|\chi_E(x)\int_F(b(x)-b(y))\frac{\Omega(x-y)}{|x-y|^n}dy\Big|\\
&=\int_{\mathbb{R}^n}\chi_{E\times F}(x,y)|b(x)-b(y)|\frac{|\Omega(x-y)|}{|x-y|^n}dy\\
&\gtrsim w_{2^{-n-2}}(b;Q)|Q|^{-1}\int_{\mathbb{R}^n}\chi_{G}(x,y)|\Omega(x-y)|dy\\
&\geq\epsilon_0w_{2^{-n-2}}(b;Q)|Q|^{-1}|G_x|,
\end{align*}
where $G_x:=\{y:(x,y)\in G\}$. From the proof of Lemma \ref{mean value method} in \cite{LeOmRi}, it is not hard to check that
\begin{align}\label{measure G and Q}
|G_x|\gtrsim|Q|.
\end{align}
Hence,
$$\chi_E(x)|[b,T_{\Omega}](\chi_F)(x)|\geq \mu w_{2^{-n-2}}(b;Q)$$
for $\mu$ small enough.

From the assumption, for any $t>0$,
\begin{align*}
|\{x\in\mathbb{R}^n:|[b,T_\Omega](\chi_F)(x)|>t|
\lesssim\int_F\Phi\Big(\frac{\alpha}{t}\Big)dx.
\end{align*}
Take $t=\mu w_{2^{-n-2}}(b;Q)/2$, we obatin
\begin{align}\label{vital1}
|E|\lesssim\frac{\alpha}{t}\log\Big(e+\frac{\alpha}{t}\Big)|F|.
\end{align}
Note that $\lim_{t\rightarrow\infty}\log(e+1/t)/t=0$, thanks to \eqref{vital1} and $|E|\sim|F|$, we get
\begin{align*}
t=\mu w_{2^{-n-2}}(b;Q)/2\lesssim1.
\end{align*}
Thus, the conclusion follows by \eqref{transform to local mean oscillation}.
\end{proof}

\begin{proof}[Proof of Theorem \ref{characterization endpoint estimate of commutator of VCZO with b real}]
It suffices to prove $(2)\Rightarrow(1)$ since $(1)\Rightarrow(2)$ is proved in \cite{WenWuXue}. Given a cube $Q$, let $E$, $F$, $G$ be the sets stated in Lemma \ref{mean value method}. Note that the supremum of $\mathcal{V}_\rho$ is taken over all sequence $\{\varepsilon_i\}$ decreasing to zero, this leads to our being able to choose $\{\varepsilon_i\}$ to satisfy $$\sharp\{i\in\mathbb{N}:\{y:\varepsilon_{i+1}<|x-y|\leq\varepsilon_i\}\cap F\neq\emptyset\}\lesssim 1$$
and
$$\Big(\bigcup_{i=1}^{\infty}\{y:\varepsilon_{i+1}<|x-y|\leq\varepsilon_i\}\Big)\bigcap F=F,$$
where $\sharp$ is the counting measure. So there are only finite non-zero terms in the following summation, Lemma \ref{mean value method} entails that
\begin{align}\label{use for variation}
&\Big(\sum_{i=1}^{\infty}\Big|\int_{\varepsilon_{i+1}<|x-y|\leq\varepsilon_i}\chi_E(x)(b(x)-b(y))
\Omega(x-y)|x-y|^{-n}
\chi_F(y)dy\Big|^\rho\Big)^{1/\rho}\\
&\quad\sim\sum_{i=1}^{\infty}\Big|\int_{\varepsilon_{i+1}<|x-y|\leq\varepsilon_i}\chi_E(x)(b(x)-b(y))
\Omega(x-y)|x-y|^{-n}
\chi_F(y)dy\Big|\nonumber\\
&\quad=\int_{F}|b(x)-b(y)||\Omega(x-y)||x-y|^{-n}\chi_E(x)dy.\nonumber
\end{align}

By \eqref{measure G and Q}, \eqref{use for variation} and Lemma \ref{mean value method},
\begin{align*}
&\mathcal{V}_\rho(\mathcal{T}_{\Omega,b})(\chi_F)(x)\chi_E(x)\\
&\quad=
\sup_{\varepsilon_i\downarrow0}\Big(\sum_{i=1}^\infty\Big|
\int_{\varepsilon_{i+1}<|x-y|\leq\varepsilon_i}(b(x)-b(y))\Omega(x-y)|x-y|^{-n}\chi_E(x)
\chi_F(y)dy\Big|^\rho\Big)^{1/\rho}\\
&\quad\gtrsim
\int_{\mathbb{R}^n}\chi_G(x,y)|b(x)-b(y)||\Omega(x-y)||x-y|^{-n}dy\\
&\quad\gtrsim\varepsilon_0 w_{2^{-n-2}}(b;Q)|G_x||Q|^{-1}\gtrsim\varepsilon_0 w_{2^{-n-2}}(b;Q).
\end{align*}
Hence, $\chi_E(x)|\mathcal{V}_\rho(\mathcal{T}_{\Omega,b})(\chi_F)(x)|\geq \mu w_{2^{-n-2}}(b;Q)$ for $\mu$ small enough.

From the assumption, for any $t>0$,
\begin{align*}
\omega\Big(\Big\{x\in\mathbb{R}^n:|\mathcal{V}_\rho(\mathcal{T}_{\Omega,b})(\chi_F)(x)|>t\Big\}\Big)
\lesssim\int_F\Phi\Big(\frac{1}{t}\Big)\omega(x)dx.
\end{align*}
Take $t=\mu w_{2^{-n-2}}(b;Q)/2$, we obatin
\begin{align}\label{vital2}
w(E)\lesssim\frac{1}{t}\log\Big(e+\frac{1}{t}\Big)w(F).
\end{align}
Observe that $\omega(k_0Q)\sim\omega(Q)$ by the doubling property of $w$, combining the fact
\begin{align}\label{important weight property}
\frac{|E|}{|Q|}\leq[w]_{A_1}\frac{w(E)}{w(Q)},
\end{align}
see \cite[Lemma 7.2.1]{Gra}, we infer that
\begin{align}\label{vital3}
\frac{\omega(E)}{\omega(F)}\geq\frac{\omega(E)}{\omega(k_0Q)}\sim\frac{\omega(E)}{\omega(Q)}
\geq\frac{|E|}{[\omega]_{A_1}|Q|}\sim\frac{1}{[\omega]_{A_1}},
\end{align}
where we use $|E|\sim|Q|$. By making use of $\lim_{t\rightarrow\infty}\log(e+1/t)/t=0$, \eqref{vital2} and \eqref{vital3}, we achieve
\begin{align*}
t=\mu w_{2^{-n-2}}(b;Q)/2\lesssim1.
\end{align*}
Thus, the conclusion follows by \eqref{transform to local mean oscillation}.
\end{proof}

\begin{proof}[Proof of Theorem \ref{characterization endpoint estimate of commutator of V approximate}]
$(1)\Rightarrow(2)$ is a consquence of results in \cite{GuoWenWuYang,LeOmRi2017}. We concentrate on proving $(2)\Rightarrow(1)$. Without loss of generality, we assume that $b$ and $\phi$ are real-valued,
and $\phi(z)\geq 1$ for $z\in B(z_0,\delta)$, where $|z_0|=1$ and $\delta>0$ is a small constant. For arbitrary cube $Q$, we denote
\begin{align*}
P:=Q-10\sqrt{n}\delta^{-1}l_Qz_0.
\end{align*}
According to the definition of $w_{2^{-n-2}}(b;Q)$, there is a set $\tilde Q\subset Q$, such that $|\tilde Q|=2^{-n-2}|Q|$ and
\begin{align*}
w_{2^{-n-2}}(b;Q)\leq |b(x)-m_b(P)|,~x\in \widetilde Q.
\end{align*}
Then, ultimate the definition of $m_b(P)$, there are subsets $E\subset \tilde Q$ and $F\subset P$ satisfying
\begin{align}\label{add1}
  |E|=|\tilde Q|/2=|Q|/2^{n+3},\ |F|=|P|/2=|Q|/2,
\end{align}
  and
\begin{align}\label{add2}
  w_{2^{-n-2}}(b;Q)\leq |b(x)-b(y)|,\quad\forall\, x\in E,\, y\in F,
\end{align}
and $b(x)-b(y)$ does not change sign in $E\times F$.

For $x\in E\subset Q, y\in F\subset P$, we claim that
\begin{align}\label{function property}
\phi_{10\sqrt{n}\delta^{-1}l_Q}(x-y)\gtrsim \frac{1}{|Q|},\,\lim_{t\rightarrow 0}\left|\phi_{t}(x-y))\right|=0.
\end{align}
In fact,
\begin{align*}
x-y\in 2l_QQ_0+10\sqrt{n}\delta^{-1}l_Qz_0\subset (l_QQ_0)^c,\ \ \ \frac{x-y}{10\sqrt{n}\delta^{-1}l_Q}\in \frac{\delta}{5\sqrt{n}}Q_0+z_0\subset B(z_0,\delta),
\end{align*}
where $Q_0:[-1/2,1/2]^n$. Recall that $\phi(z)\geq1$ for $z\in B(z_0,\delta)$. Thus, for $x\in E$, $y\in F$, we end up with
\begin{align*}
  \phi_{10\sqrt{n}\delta^{-1}l_Q}(x-y)
  \gtrsim \frac{1}{|Q|}\phi\Big(\frac{x-y}{10\sqrt{n}\delta^{-1}l_Q}\Big)\geq\frac{1}{|Q|},
\end{align*}
and
\begin{align*}
  \lim_{t\rightarrow 0}\left|\phi_{t}(x-y))\right|
  &=
  \lim_{t\rightarrow 0}\frac{1}{t^n}\left|\phi\Big(\frac{x-y}{t}\Big)\right|\lesssim \lim_{t\rightarrow 0}\frac{1}{t^n}\Big(\frac{|x-y|}{t}\Big)^{-n-1}\\
  &\lesssim \lim_{t\rightarrow 0}t(|x-y|)^{-n-1}=0.
\end{align*}

Let $f=\chi_F$. Then for $x\in E$, by Fatou's lemma, \eqref{add1}, \eqref{add2} and \eqref{function property},
\begin{align*}
\mathcal{V}_\rho((\Phi\star f)_b)(x)&=\sup_{t_i\downarrow0}\Big(\sum_{i=1}^\infty\Big|\int_F(b(x)-b(y))(\phi_{t_i}(x-y)-\phi_{t_{i+1}}(x-y))
\Big|^\rho\Big)^{1/\rho}\\
&\geq\liminf_{t\rightarrow0}\Big|\int_F(b(x)-b(y))(\phi_{10\sqrt{n}\delta^{-1}l_Q}(x-y)-\phi_t(x-y))dy\Big|\\
&\geq\int_F|b(x)-b(y)|\liminf_{t\rightarrow0}|\phi_{10\sqrt{n}\delta^{-1}l_Q}(x-y)-\phi_t(x-y)|dy\\
&=\int_F|b(x)-b(y)||\phi_{10\sqrt{n}\delta^{-1}l_Q}(x-y)|dy\\
&\gtrsim w_{2^{-n-2}}(b;Q).
\end{align*}
The rest of the proof is the same as $(2)\Rightarrow(1)$ of Theorem \ref{characterization endpoint estimate of commutator of VCZO with b real}, except for replacing $\mathcal{V}_\rho(\mathcal{T}_{\Omega,b})(f)$ with $\mathcal{V}_\rho((\Phi\star f)_b)$.
\end{proof}

\begin{proof}[Proof of Theorem \ref{characterization endpoint estimate of commutator of Marcinkiewicz}]
$(1)\Rightarrow(2)$ is presented in \cite{Zhang}. We devoted to proving $(2)\Rightarrow(1)$. Given a cube $Q$, let $E$, $F$, $k_0$ be the sets stated in Lemma \ref{mean value method}. It is direct that
\begin{align*}
u_{\Omega,b}(\chi_F)(x)&=\Big(\int_0^\infty\Big|\int_{\{y:|x-y|\leq t\}\cap F}(b(x)-b(y))\frac{\Omega(x-y)}{|x-y|^{n-1}}dy\Big|^2\frac{dt}{t^3}\Big)^{1/2}\\
&\geq\Big(\int_{(k_0+1)\sqrt{n}l_Q}^\infty\Big|\int_{\{y:|x-y|\leq t\}\cap F}(b(x)-b(y))\frac{\Omega(x-y)}{|x-y|^{n-1}}dy\Big|^2\frac{dt}{t^3}\Big)^{1/2}
\end{align*}
For $x\in E$, observe that $\{y:|x-y|\leq t\}\cap F=F$ when $t\geq (k_0+1)\sqrt{n}l_Q$. Then using Lemma \ref{mean value method}, we have
\begin{align*}
&u_{\Omega,b}(\chi_F)(x)\chi_E(x)\\
&\quad\geq\Big(\int_{(k_0+1)\sqrt{n}l_Q}^\infty\Big|\int_{ F}(b(x)-b(y))\Omega(x-y)|x-y|^{1-n}\chi_E(x)dy\Big|^2\frac{dt}{t^{2n+1}}\Big)^{1/2}\\
&\quad\gtrsim\epsilon_0 w_{2^{-n-2}}(b;Q)|Q|^{1/n}\Big(\int_{(k_0+1)\sqrt{n}l_Q}^\infty\frac{dt}{t^{2n+1}}\Big)^{1/2}\\
&\quad\sim\epsilon_0 w_{2^{-n-2}}(b;Q).
\end{align*}
Therefore
$$\chi_E(x)u_{\Omega,b}f(x)\geq\mu w_{2^{-n-2}}(b;Q)$$
for $\mu$ small enough. The rest of the proof is the same as $(2)\Rightarrow(1)$ of Theorem \ref{characterization endpoint estimate of commutator of VCZO with b real}.
\end{proof}

\begin{proof}[Proof of Theorem \ref{characterization endpoint estimate of commutator of parameter with b real}]
$(1)\Rightarrow(2)$ is given in \cite{DiXu}. We shall prove $(2)\Rightarrow(1)$. We first consider $u_{S,b}^\varrho$. Given a cube $Q$, let $E$, $F$, $k_0$ be the sets stated in Lemma \ref{mean value method}. Observe that
\begin{align}\label{observe1}
|x-y|\leq(k_0+1)\sqrt{n}l_Q,~x,y\in E,
\end{align}
\begin{align}\label{observe2}
\{z:|y-z|\leq t\}\cap F=F,~y\in E,~t\geq (k_0+1)\sqrt{n}l_Q.
\end{align}
Let $f=\chi_F$. For $x\in E$, Lemma \ref{mean value method} implies that
\begin{align*}
u_{S,b}^\varrho (\chi_F)(x)&=\Big(\iint_{\Gamma(x)}t^{-2\varrho}\Big|\int_{|y-z|<t}\frac{\Omega(y-z)}{|y-z|^{n-\varrho}}(b(x)-b(z))
\chi_F(z)dz\Big|^2\frac{dydt}{t^{n+1}}\Big)^{1/2}\\
&\geq\Big(\int_{(k_0+1)\sqrt{n}l_Q}^\infty\int_E t^{-2\varrho}\Big|\int_{F}
\frac{\Omega(y-z)}{|y-z|^{n-\varrho}}(b(x)-b(z))\chi_F(z)dz\Big|^2\frac{dydt}{t^{n+1}}\Big)^{1/2}\\
&\gtrsim \epsilon_0w_{2^{-n-2}}(b;Q)|Q|^{\varrho/n-1}|Q|\Big(\int_{(k_0+1)\sqrt{n}l_Q}^\infty\int_E t^{-2\varrho-n-1}dydt\Big)^{1/2}\\
&\sim\epsilon_0w_{2^{-n-2}}(b;Q)|Q|^{\varrho/n}|Q|^{-\varrho/n}=\epsilon_0w_{2^{-n-2}}(b;Q).
\end{align*}
Therefore
$$\chi_E(x)u_{S,b}^\varrho (\chi_F)(x)\geq\mu w_{2^{-n-2}}(b;Q)$$
for $\mu$ small enough. The rest of the proof is the same as $(2)\Rightarrow(1)$ of Theorem \ref{characterization endpoint estimate of commutator of VCZO with b real}.

Next, we deal with $g_{\lambda,b}^{\ast,\varrho}$. We still use the notations above. Again using \eqref{observe1} and \eqref{observe2}, then for $x\in E$,
\begin{align*}
&g_{\lambda,b}^{\ast,\varrho}(\chi_F)(x)\\
&=\Big(\iint_{\mathbb{R}_+^{n+1}}\Big(\frac{t}{t+|x-y|}\Big)^{\lambda n}\Big|\frac{1}{t^{\varrho}}\int_{|y-z|<t}\frac{\Omega(y-z)}{|y-z|^{n-\varrho}}(b(x)-b(z))
\chi_F(z)dz\Big|^2\frac{dydt}{t^{n+1}}\Big)^{1/2}\\
&\geq\Big(\int_{(k_0+1)\sqrt{n}l_Q}^\infty\int_E \Big(\frac{t}{t+|x-y|}\Big)^{\lambda n}t^{-2\varrho}\Big|\int_{F}
\frac{\Omega(y-z)}{|y-z|^{n-\varrho}}(b(x)-b(z))dz\Big|^2\frac{dydt}{t^{n+1}}\Big)^{1/2}\\
&\geq2^{-\lambda n}\Big(\int_{(k_0+1)\sqrt{n}l_Q}^\infty\int_E t^{-2\varrho}\Big|\int_{F}
\frac{\Omega(y-z)}{|y-z|^{n-\varrho}}(b(x)-b(z))dz\Big|^2\frac{dydt}{t^{n+1}}\Big)^{1/2}.
\end{align*}
Owing to Lemma \ref{mean value method}, we further have
\begin{align*}
&g_{\lambda,b}^{\ast,\varrho}(\chi_F)(x)\chi_E(x)\\
&\quad\gtrsim \epsilon_0w_{2^{-n-2}}(b;Q)|Q|^{\varrho/n-1}|Q|\Big(\int_{(k_0+1)\sqrt{n}l_Q}^\infty\int_E t^{-2\varrho-n-1}dydt\Big)^{1/2}\\
&\quad\sim\epsilon_0w_{2^{-n-2}}(b;Q)|Q|^{\varrho/n}|Q|^{-\varrho/n}=\epsilon_0w_{2^{-n-2}}(b;Q).
\end{align*}
Therefore
$$\chi_E(x)g_{\lambda,b}^{\ast,\varrho}(\chi_F)(x)\geq\mu w_{2^{-n-2}}(b;Q)$$
for $\mu$ small enough. The rest of the proof is the same as $(2)\Rightarrow(1)$ of Theorem \ref{characterization endpoint estimate of commutator of VCZO with b real}.
\end{proof}
\section{Characterization of the weighted endpoint estimate of commutator of non-degenerate Calder\'{o}n-Zygmund operator}
In this section, we prove Theorem \ref{necesity endpoint estimate of commutator of non-degenerate CZO}. We avoid using the mean value method so that we can relax the restriction that the pointwise multiplier be real-valued.

\begin{lemma}{\rm(\cite{Hy})}\label{necessity complex}
Let $K$ be a non-degenerate Calder\'{o}n-Zygmund kernel. Then for every $A\geq3$ and every ball $B:=B(y_0,r)$, there is a disjoint ball $\widetilde{B}:=B(x_0,r)$ at distance dist$(B,\widetilde{B})\sim Ar$ such that $|K(x_0,y_0)|\sim (Ar)^{-n}$, and for all $y_1\in B$ and $x_1\in\widetilde{B}$, we have
\begin{align*}
\int_{B}|K(x_1,y)-K(x_0,y_0)|dy+\int_{\widetilde{B}}|K(x,y_1)-K(x_0,y_0)|dx\lesssim\epsilon_AA^{-n},
\end{align*}
where $\epsilon_A\rightarrow0$ as $A\rightarrow\infty$ and implicit constants depend at most on $c_K$, $\omega$ and $n$, as well as $c_0$ or $\Omega(\theta_0)$ from Definition \ref{non-degenerate CZ}.
\end{lemma}

\begin{proof}[Proof of Theorem \ref{necesity endpoint estimate of commutator of non-degenerate CZO}]
Fix a ball $B:=B(y_0,r)$. Based on Lemma \ref{necessity complex} and the doubling property of $w$, there is a disjoint ball $\widetilde{B}:=B(x_0,r)$ satisfying
\begin{align}\label{weight B compare}
w(B)\lesssim w((A+1)\widetilde{B})\leq C_ww(\widetilde{B}),
\end{align}
\begin{align}\label{size of special kernel}
|K(x_0,y_0)|\sim (Ar)^{-n},
\end{align}
\begin{align}\label{uniformly smooth kernel1}
\sup_{x\in\widetilde{B}}\int_{B}|K(x,y)-K(x_0,y_0)|dy\lesssim\epsilon_AA^{-n},
\end{align}
\begin{align}\label{uniformly smooth kernel2}
\sup_{y\in B}\int_{\widetilde{B}}|K(x,y)-K(x_0,y_0)|dx\lesssim\epsilon_AA^{-n},
\end{align}
where $\epsilon_A\rightarrow0$ as $A\rightarrow\infty$, and implicit constants depend at most on $c_K$, $\omega$ and $n$, as well as $c_0$ or $\Omega(\theta_0)$ from Definition \ref{non-degenerate CZ}. Hence, by Tonelli's theorem and \eqref{uniformly smooth kernel2},
\begin{align}\label{use for estimate measure}
&\int_{\widetilde{B}}\int_B|b(y)-b_B||K(x,y)-K(x_0,y_0)|dydx\\
&\quad=\int_B|b(y)-b_B|\int_{\widetilde{B}}|K(x,y)-K(x_0,y_0)|dxdy\nonumber\\
&\quad\lesssim\epsilon_AA^{-n}\int_B|b(y)-b_B|dy\nonumber\\
&\quad=\epsilon_AA^{-n}|B|(|b(y)-b_B|)_B.\nonumber
\end{align}
Denote
\begin{align*}
D:=\Big\{x\in\widetilde{B}:\int_B|b(y)-b_B||K(x,y)-K(x_0,y_0)|dy\leq4c\epsilon_AA^{-n}(|b(y)-b_B|)_B\Big\}.
\end{align*}
Then \eqref{use for estimate measure} implies that $|\widetilde{B}\backslash D|\leq|B|/4=|\widetilde{B}|/4$, which yields that
\begin{align}\label{use for find point}
w(D)\geq c_ww(\widetilde{B})
\end{align}
by \eqref{important weight property}.

Define $f_1=\chi_B$ and $f_2=(g-g_B)\chi_B$, where
\[
g(y)\\
:=
\begin{cases}
\frac{\overline{b(y)-b_B}}{|b(y)-b_B|},& b(y)\neq b_B,\\
0,& b(y)=b_B.
\end{cases}
\]
By the assumption,
\begin{align*}
w(\{x\in\mathbb{R}^n:|[b,T]f_1(x)|>\lambda\alpha\})\leq C_w\int_B\Phi\Big(\frac{1}{\lambda}\Big)w(x)dx=C_w\Phi\Big(\frac{1}{\lambda}\Big)w(B)
\end{align*}
and
\begin{align*}
w(\{x\in\mathbb{R}^n:|[b,T]f_2(x)|>\lambda\alpha\})\leq C_w\int_B\Phi\Big(\frac{2}{\lambda}\Big)w(x)dx=C_w\Phi\Big(\frac{2}{\lambda}\Big)w(B),
\end{align*}
where we use $\|f_2\|_{L^\infty}\leq2$ in the second inequality. Then by making use of \eqref{weight B compare}, there exists $\lambda_0>0$ big enough such that
\begin{align*}
&w(\{x\in\mathbb{R}^n:|[b,T]f_1(x)|>\lambda_0\alpha\})+w(\{x\in\mathbb{R}^n:|[b,T]f_2(x)|>\lambda_0\alpha\})\\
&\quad\leq C_ww(B)\Big[\Phi\Big(\frac{1}{\lambda_0}\Big)+\Phi\Big(\frac{2}{\lambda_0}\Big)\Big]\leq \frac{c_w}{2}w(\widetilde{B}).
\end{align*}
This, together with \eqref{use for find point}, shows that
\begin{align*}
&w\big(D\backslash(\{x\in\mathbb{R}^n:|[b,T]f_1(x)|>\lambda_0\alpha\}\cup
\{x\in\mathbb{R}^n:|[b,T]f_2(x)|>\lambda_0\alpha\})\big)\\
&\quad\geq w(D)-w(\{x\in\mathbb{R}^n:|[b,T]f_1(x)|>\lambda_0\alpha\})-w(\{x\in\mathbb{R}^n:|[b,T]f_2(x)|>\lambda_0\alpha\})\\
&\quad\geq\frac{c_w}{2}w(\widetilde{B}).
\end{align*}
Therefore, there exists $\widetilde{x}\in\widetilde{B}$ such that
\begin{align}\label{point1}
\int_B|b(y)-b_B||K(\widetilde{x},y)-K(x_0,y_0)|dy\leq4c\epsilon_AA^{-n}(|b-b_B|)_B,
\end{align}
\begin{align}\label{point2}
\max\{|[b,T]f_1(\widetilde{x})|,|[b,T]f_2(\widetilde{x})|\}\leq\lambda_0\alpha.
\end{align}

Now for $\widetilde{x}\in\widetilde{B}$, we write
\begin{align*}
[b,T]f_1(\widetilde{x})&=\int_B(b(\widetilde{x})-b(y))(K(\widetilde{x},y)-K(x_0,y_0)+K(x_0,y_0))dy\\
&=K(x_0,y_0)|B|(b(\widetilde{x})-b_B)+\int_B(b(\widetilde{x})-b(y))(K(\widetilde{x},y)-K(x_0,y_0))dy.
\end{align*}
By \eqref{size of special kernel}, \eqref{uniformly smooth kernel1}, \eqref{point1} and \eqref{point2}, for $\widetilde{x}\in\widetilde{B}$, we deduce that
\begin{align*}
\lambda_0\alpha&\geq|[b,T]f_1(\widetilde{x})|\\
&\geq|K(x_0,y_0)||B||b(\widetilde{x})-b_B|-\int_B|b(\widetilde{x})-b(y)||K(\widetilde{x},y)-K(x_0,y_0)|dy\\
&\geq|K(x_0,y_0)||B||b(\widetilde{x})-b_B|-|b(\widetilde{x})-b_B|\int_B|K(\widetilde{x},y)-K(x_0,y_0)|dy\\
&\quad-\int_B|b(y)-b_B||K(\widetilde{x},y)-K(x_0,y_0)|dy\\
&\gtrsim A^{-n}|b(\widetilde{x})-b_B|-\epsilon_AA^{-n}|b(\widetilde{x})-b_B|-4c\epsilon_AA^{-n}(|b-b_B|)_B.
\end{align*}
Rearranging the terms on both sides of the inequality yields
\begin{align}\label{final estimate1}
|b(\widetilde{x})-b_B|\lesssim \frac{A^n\lambda_0\alpha}{1-\epsilon_A}+\frac{4c\epsilon_A}{1-\epsilon_A}(|b-b_B|)_B,
\end{align}
where the implicit constant is independent of $B$.

On the other hand, note that
\begin{align*}
\int_Bf_2(y)dy=0,
\end{align*}
then
\begin{align*}
\int_Bb(y)f_2(y)dy=\int_B(b(y)-b_B)f_2(y)dy=\int_B|b(y)-b_B|dy.
\end{align*}
From this, we can write
\begin{align*}
[b,T]f_2(\widetilde{x})&=\int_B(b(\widetilde{x})-b(y))f_2(y)K(\widetilde{x},y)dy\\
&=\int_B(b(\widetilde{x})-b(y))f_2(y)(K(\widetilde{x},y)-K(x_0,y_0)+K(x_0,y_0))dy\\
&=\int_B(b(\widetilde{x})-b(y))f_2(y)(K(\widetilde{x},y)-K(x_0,y_0))dy\\
&\quad+\int_B(b(\widetilde{x})-b(y))f_2(y)K(x_0,y_0)dy\\
&=\int_B(b(\widetilde{x})-b(y))f_2(y)(K(\widetilde{x},y)-K(x_0,y_0))dy\\
&\quad+b(\widetilde{x})K(x_0,y_0)\int_Bf_2(y)dy-K(x_0,y_0)\int_Bb(y)f_2(y)dy\\
&=\int_B(b(\widetilde{x})-b(y))f_2(y)(K(\widetilde{x},y)-K(x_0,y_0))dy\\
&\quad-K(x_0,y_0)\int_B|b(y)-b_B|dy.
\end{align*}
In virtue of \eqref{uniformly smooth kernel1}, \eqref{point1} and $\|f_2\|_{L^\infty}\leq2$, we get
\begin{align*}
&\Big|\int_B(b(\widetilde{x})-b(y))f_2(y)(K(\widetilde{x},y)-K(x_0,y_0))dy\Big|\\
&\quad\leq2\int_B|b(\widetilde{x})-b(y)||K(\widetilde{x},y)-K(x_0,y_0)|dy\\
&\quad\leq2\int_B|b(\widetilde{x})-b_B||K(\widetilde{x},y)-K(x_0,y_0)|dy\\
&\qquad+2\int_B|b(y)-b_B||K(\widetilde{x},y)-K(x_0,y_0)|dy\\
&\quad\lesssim2\epsilon_AA^{-n}|b(\widetilde{x})-b_B|+8c\epsilon_AA^{-n}(|b-b_B|)_B.
\end{align*}
Then by \eqref{point2}, we infer that
\begin{align*}
\lambda_0\alpha&\geq|[b,T]f_2(\widetilde{x})|\\
&\gtrsim A^{-n}r^{-n}\int_B|b(y)-b_B|dy-2\epsilon_AA^{-n}|b(\widetilde{x})-b_B|-8c\epsilon_AA^{-n}(|b-b_B|)_B\\
&\sim(A^{-n}-8c\epsilon_AA^{-n})(|b-b_B|)_B-2\epsilon_AA^{-n}|b(\widetilde{x})-b_B|,
\end{align*}
which yields that
\begin{align*}
\frac{1}{|B|}\int_B|b(y)-b_B|dy\lesssim A^n\lambda_0\alpha+2\epsilon_A|b(\widetilde{x})-b_B|+8c\epsilon_A
\frac{1}{|B|}\int_B|b(y)-b_B|dy.
\end{align*}
Plugging \eqref{final estimate1} into the above estimate gives
\begin{align*}
\frac{1}{|B|}\int_B|b(y)-b_B|dy\leq C_1\Big( A^n\lambda_0\alpha+\frac{2\epsilon_AA^n\lambda_0\alpha}{1-\epsilon_A}+\frac{8c\epsilon_A}{1-\epsilon_A}
\frac{1}{|B|}\int_B|b(y)-b_B|dy\Big),
\end{align*}
where $C_1,c$ is independent of $B$. Since $\epsilon_A\rightarrow0$ as $A\rightarrow\infty$, we can take $A$ large enough such that $1>\epsilon_A+16cC_1\epsilon_A$. Note that $b\in L_{loc}^1(\mathbb{R}^n)$, we deduce that
\begin{align*}
\frac{1}{|B|}\int_B|b(y)-b_B|dy\leq2C_1A^n\lambda_0\alpha+\frac{4C_1\epsilon_AA^n\lambda_0\alpha}{1-\epsilon_A}.
\end{align*}
This accomplish the proof of Theorem \ref{necesity endpoint estimate of commutator of non-degenerate CZO} by taking the supremum over all $B\subset\mathbb{R}^n$.
\end{proof}

\section{The weighted endpoint estimate of commutator of fractional integral does not imply $b\in {\rm BMO}$}
In this section, we devoted to proving Theorem \ref{counterexample}. We first prove some necessary lemmas.

\begin{lemma}\label{use for key lemma1}
Let $0<\alpha<n$ and $B$ be a measurable set. Then for any $x\in\mathbb{R}^n$, there exists $C_{n,\alpha}>0$ such that
\begin{align*}
I_\alpha(\chi_B)(x)\leq C_{n,\alpha}|B|^{\alpha/n}.
\end{align*}
\end{lemma}
\begin{proof}
We decompose
\begin{align*}
I_\alpha(\chi_B)(x)&=\int_{B\cap\{y:|x-y|\geq|B|^{1/n}\}}|x-y|^{\alpha-n}dy+
\int_{B\cap\{y:|x-y|<|B|^{1/n}\}}|x-y|^{\alpha-n}dy\\
&=:I_\alpha^{global}(\chi_B)(x)+I_\alpha^{local}(\chi_B)(x).
\end{align*}
For $I_\alpha^{global}(\chi_B)$, it is direct that
\begin{align*}
I_\alpha^{global}(\chi_B)(x)\leq|B|^{(\alpha-n)/n}|B|=|B|^{\alpha/n}.
\end{align*}
For $I_\alpha^{local}(\chi_B)$, polar coordinates gives
\begin{align*}
I_\alpha^{local}(\chi_B)(x)&\leq\int_{|x-y|<|B|^{1/n}}|x-y|^{\alpha-n}dy\\
&=c_n\int_{0}^{|B|^{1/n}}r^{\alpha-n}r^{n-1}dr=c_{n,\alpha}|B|^{\alpha/n}.
\end{align*}
Combined with the estimates of $I_\alpha^{global}(\chi_B)$ and $I_\alpha^{local}(\chi_B)$, we get the desired result.
\end{proof}

\begin{lemma}\label{key lemma1}
Let $\{Q_j\}_{j\in\mathbb{N}}$ be a collection of cubes with $l(Q_j)=e^{-4^j}$. Let $0<\alpha<n$, $q=n/(n-\alpha)$, $\Phi(t)=t\log(e+t)$ and $w\in A_1$. Then there exists $C_w>0$ such that
\begin{align*}
w(\{x\in\mathbb{R}^n:jI_\alpha(\chi_{Q_j}|f|)(x)>1/2\})\leq C_w\Big(j4^{-j}\int_{\mathbb{R}^n}\Phi(|f(x)|)w(x)^{1/q}dx\Big)^q.
\end{align*}
\end{lemma}
\begin{proof}
Let $C_{n,\alpha}$ be as in Lemma \ref{use for key lemma1}. We decompose $\chi_{Q_j}|f|$ by
\begin{align*}
\chi_{Q_j}(y)|f(y)|&=\chi_{Q_j}(y)|f(y)|
\chi_{\{y\in\mathbb{R}^n:|f(y)|\leq(4jC_{n,\alpha}|Q_j|^{\alpha/n})^{-1}\}}(y)\\
&\quad+\chi_{Q_j}(y)|f(y)|
\chi_{\{y\in\mathbb{R}^n:|f(y)|>(4jC_{n,\alpha}|Q_j|^{\alpha/n})^{-1}\}}(y)\\
&=:f_j^{local}(y)+f_j^{global}(y).
\end{align*}
By Lemma \ref{use for key lemma1},
\begin{align*}
jI_\alpha(f_j^{local})&\leq(4jC_{n,\alpha}|Q_j|^{\alpha/n})^{-1}jI_\alpha(\chi_{Q_j})\\
&\leq(4jC_{n,\alpha}|Q_j|^{\alpha/n})^{-1}jC_{n,\alpha}|Q_j|^{\alpha/n}=1/4.
\end{align*}
Then
\begin{align*}
&w(\{x\in\mathbb{R}^n:jI_\alpha(\chi_{Q_j}|f|)(x)>1/2\})\\
&\quad\leq w(\{x\in\mathbb{R}^n:jI_\alpha(f_j^{local})(x)>1/4\})+
w(\{x\in\mathbb{R}^n:jI_\alpha(f_j^{global})(x)>1/4\})\\
&\quad=w(\{x\in\mathbb{R}^n:I_\alpha(f_j^{global})(x)>1/(4j)\}).
\end{align*}
Hence, it remains to prove the desired result for
\begin{align}\label{final estimate of key lemma1}
w(\{x\in\mathbb{R}^n:I_\alpha(f_j^{global})(x)>1/(4j)\}).
\end{align}
Before this, we first claim that for all $j\geq1$,
\begin{align}\label{log estimate1}
\log(e+(4jC_{n,\alpha}|Q_j|^{\alpha/n})^{-1})\geq c_\alpha4^j.
\end{align}
In fact, observe that there exists $j_0\geq1$ such that for every $j>j_0$, we have
\begin{align*}
\log(e+(4jC_{n,\alpha}|Q_j|^{\alpha/n})^{-1})&\geq\log(4jC_{n,\alpha}|Q_j|^{\alpha/n})^{-1}\\
&=\log(4jC_{n,\alpha}e^{-\alpha4^j})^{-1}\\
&=\alpha4^j-\log j-\log 4C_{n,\alpha}\geq\alpha4^j/2.
\end{align*}
while for $1\leq j\leq j_0$,
\begin{align*}
\log(e+(4jC_{n,\alpha}|Q_j|^{\alpha/n})^{-1})\geq1\geq4^{-j_0}4^j.
\end{align*}
Thus, \eqref{log estimate1} holds for any $j\geq1$.

We return to consider \eqref{final estimate of key lemma1}. According to the weighted endpoint estimate of $I_\alpha$ due to Muckenhoupt and  Wheeden \cite{MuWh} and \eqref{log estimate1},
\begin{align*}
&w(\{x\in\mathbb{R}^n:I_\alpha(f_j^{global})(x)>1/(4j)\})\\
&\quad\leq C_w\Big(4j\int_{\mathbb{R}^n}f_j^{global}(x)w(x)^{1/q}dx\Big)^q\\
&\quad=C_w\Big(4j\int_{\mathbb{R}^n}\frac{f_j^{global}(x)\log(e+f_j^{global}(x))}
{\log(e+f_j^{global}(x))}w(x)^{1/q}dx\Big)^q\\
&\quad\leq C_w\Big(\frac{4j}{\log(e+(4jC_{n,\alpha}|Q_j|^{-\alpha/n})^{-1})}
\int_{\mathbb{R}^n}\Phi(|f(x)|)w(x)^{1/q}dx\Big)^q\\
&\quad\leq C_w\Big(j4^{-j}\int_{\mathbb{R}^n}\Phi(|f(x)|)w(x)^{1/q}dx\Big)^q.
\end{align*}
This finish the proof of Lemma \ref{key lemma1}.
\end{proof}

\begin{lemma}\label{A1 theorem}
Let $\{Q_j\}_{j\in\mathbb{N}}$ be a collection of cubes with $l(Q_j)=e^{-4^j}$, that is contained in $Q_0:=[-1/2,1/2]^n$. If $w\in A_1$, then there exist constants $C_w',C_w'',\gamma_w>0$ such that
\begin{align*}
w(x)\geq C_w'(1+|x|)^{-n},~~w(Q_j)\leq C_w''e^{-\gamma_w4^j}.
\end{align*}
\end{lemma}
\begin{proof}
Let $x\in\mathbb{R}^n$ and $B:=B(0,1+|x|)$. Let $M$ be the Hardy-Littlewood maximal operator. By $w\in A_1$,
\begin{align*}
\frac{1}{|B|}\int_{B(0,1)}w(y)dy\leq\frac{1}{|B|}\int_{B}w(y)dy\leq Mw\leq[w]_{A_1}w(x),
\end{align*}
which yields that $w(x)\geq C_w'(1+|x|)^{-n}$.

In view of \cite[Proposition 7.2.8]{Gra}, there exist constants $C_w,\gamma_w>0$ such that
\begin{align*}
\frac{w(Q_j)}{w(Q_0)}\leq C_w\Big(\frac{|Q_j|}{|Q_0|}\Big)^{\gamma_w},
\end{align*}
which gives $w(Q_j)\leq C_w''e^{-\gamma_w4^j}$.
\end{proof}

\begin{lemma}\label{key lemma2}
Let $\{Q_j\}_{j\in\mathbb{N}}$ be a collection of disjoint cubes with $l(Q_j)=e^{-4^j}$, that is contained in $Q_0:=[-1/2,1/2]^n$. Let $0<\alpha<n$, $q=n/(n-\alpha)$, $\Phi(t)=t\log(e+t)$ and $w\in A_1$. Then there exists $C_{w,n,\alpha}>0$ such that
\begin{align*}
w(\{x\in\mathbb{R}^n:j\chi_{Q_j}(x)I_\alpha(|f|)(x)>1/2\})\leq C_{w,n,\alpha}\Big(j4^{-j}\int_{\mathbb{R}^n}\Phi(|f(x)|)w(x)^{1/q}dx\Big)^q.
\end{align*}
\end{lemma}
\begin{proof}
Without loss of generality, we may assume that $0<\int_{\mathbb{R}^n}\Phi(|f|)\omega^{1/q}<\infty$. Otherwise, the result is trivial. We prove this lemma by considering two cases:
\begin{align*}
\int_{\mathbb{R}^n}\Phi(|f|)\omega^{1/q}\geq e^{-\gamma_w4^j/(2q)}
\end{align*}
and
\begin{align*}
\int_{\mathbb{R}^n}\Phi(|f|)\omega^{1/q}<e^{-\gamma_w4^j/(2q)},
\end{align*}
where $\gamma_w$ is given as in Lemma \ref{A1 theorem}.

\textbf{Case 1. $\int_{\mathbb{R}^n}\Phi(|f|)\omega^{1/q}\geq e^{-\gamma_w4^j/(2q)}$}:\quad Taking advantage of Lemma \ref{A1 theorem} and $e^{-x}\lesssim x^{-\beta}$ for any $\beta>0$, we deduce that
\begin{align*}
&w(\{x\in\mathbb{R}^n:j\chi_{Q_j}(x)I_\alpha(|f|)(x)>1/2\})\\
&\quad\leq w(Q_j)\leq C_{w}''e^{-\gamma_w4^j}\\
&\quad=C_{w}''e^{-\gamma_w4^j/2}e^{-\gamma_w4^j/2}\\
&\quad\leq C_{w}\Big(\int_{\mathbb{R}^n}\Phi(|f(x)|)w(x)^{1/q}dx\Big)^q4^{-jq}j^q\\
&\quad=C_w\Big(j4^{-j}\int_{\mathbb{R}^n}\Phi(|f(x)|)w(x)^{1/q}dx\Big)^q.
\end{align*}

\textbf{Case 2. $\int_{\mathbb{R}^n}\Phi(|f|)\omega^{1/q}<e^{-\gamma_w4^j/(2q)}$}:\quad We decompose
\begin{align*}
|f|&=|f|\chi_{\{|f|\leq\delta j^{-q}(\int_{\mathbb{R}^n}\Phi(|f|)w^{1/q})^{1-q}\}}+|f|\chi_{\{|f|>\delta j^{-q}(\int_{\mathbb{R}^n}\Phi(|f|)w^{1/q})^{1-q}\}}\\
&=:f_j^{low}+f_j^{high},
\end{align*}
where $\delta$ is determined later. Then
\begin{align*}
&w(\{x\in\mathbb{R}^n:j\chi_{Q_j}(x)I_\alpha(|f|)(x)>1/2\})\\
&\quad\leq w(\{x\in\mathbb{R}^n:j\chi_{Q_j}(x)I_\alpha(f_j^{low})(x)>1/4\})\\
&\qquad+w(\{x\in\mathbb{R}^n:j\chi_{Q_j}(x)I_\alpha(f_j^{high})(x)>1/4\}).
\end{align*}

We claim that there exists $j_0\geq1$ such that for any $j>j_0$,
\begin{align}\label{case 2 claim}
w(\{x\in\mathbb{R}^n:j\chi_{Q_j}(x)I_\alpha(f_j^{low})(x)>1/4\})=0.
\end{align}
Actually, for $x\in Q_0$,
\begin{align}\label{low part decompose of fractional}
I_\alpha(f_j^{low})(x)=\int_{4Q_0}\frac{f_j^{low}(y)}{|x-y|^{n-\alpha}}dy+
\int_{(4Q_0)^c}\frac{f_j^{low}(y)}{|x-y|^{n-\alpha}}dy.
\end{align}
For $x\in Q_0$ and $y\in(4Q_0)^c$, it is not hard to check that $|x-y|\geq c_n(1+|y|)$. From this and $n-\alpha=n/q$ Lemma \ref{A1 theorem},
\begin{align}
&\int_{(4Q_0)^c}\frac{f_j^{low}(y)}{|x-y|^{n-\alpha}}dy\\
&\quad\leq c_n\int_{(4Q_0)^c}\frac{f(y)}{(1+|y|)^{n-\alpha}}dy\nonumber\\
&\quad\leq C_{w,n,\alpha}'\int_{(4Q_0)^c}|f(y)|w(y)^{1/q}dy\nonumber\\
&\quad\leq C_{w,n,\alpha}'\int_{\mathbb{R}^n}\Phi(|f(y)|)\omega(y)^{1/q}dy<C_{w,n,\alpha}'e^{-\gamma_w4^j/(2q)}.\nonumber
\end{align}
On the other hand,
\begin{align}
&\int_{4Q_0}\frac{f_j^{low}(y)}{|x-y|^{n-\alpha}}dy\\
&\quad=\int_{4Q_0\cap\{y:|x-y|\leq\delta^{-1/n}j^{q/n}(\int_{\mathbb{R}^n}\Phi(|f|)\omega^{1/q})^{q/n}\}}
\frac{f_j^{low}(y)}{|x-y|^{n-\alpha}}dy\nonumber\\
&\qquad+\int_{4Q_0\cap\{y:|x-y|>\delta^{-1/n}j^{q/n}(\int_{\mathbb{R}^n}\Phi(|f|)\omega^{1/q})^{q/n}\}}
\frac{f_j^{low}(y)}{|x-y|^{n-\alpha}}dy\nonumber\\
&\quad=:I(x)+II(x).\nonumber
\end{align}
For $I$, one can check that
\begin{align}\label{I}
I(x)&\leq\delta j^{-q}\Big(\int_{\mathbb{R}^n}\Phi(|f|)w^{1/q}\Big)^{1-q}
\int_{|x-y|\leq\delta^{-1/n}j^{q/n}(\int_{\mathbb{R}^n}\Phi(|f|)\omega^{1/q})^{q/n}}|x-y|^{\alpha-n}dy\\
&=c_{n,\alpha}\delta j^{-q}\Big(\int_{\mathbb{R}^n}\Phi(|f|)w^{1/q}\Big)^{1-q}
\Big[\delta^{-1/n}j^{q/n}\Big(\int_{\mathbb{R}^n}\Phi(|f|)\omega^{1/q}\Big)^{q/n}\Big]^{\alpha}\nonumber\\
&=c_{n,\alpha}\delta^{1/q} j^{-1}.\nonumber
\end{align}
Note that for $y\in 4Q_0$, Lemma \ref{A1 theorem} implies that $w(y)\geq c_w'$. Then
\begin{align*}
II(x)&\leq \Big[\delta^{-1/n}j^{q/n}\Big(\int_{\mathbb{R}^n}\Phi(|f|)w^{1/q}\Big)^{q/n}\Big]^{\alpha-n}
\int_{\mathbb{R}^n}\Phi(|f(y)|)w(y)^{1/q}w(y)^{-1/q}dy\\
&=c_{w,n,\alpha}'\delta^{1/q} j^{-1}.
\end{align*}
This, together with \eqref{low part decompose of fractional}-\eqref{I}, we conclude that for any $x\in Q_0$,
\begin{align}\label{estimate of low part of fractional}
I_\alpha(f_j^{low})(x)\leq C_{w,n,\alpha}'(\delta^{1/q}j^{-1}+e^{-\gamma_w4^j/(2q)}).
\end{align}
Now, we are in the position to prove \eqref{case 2 claim}. Choose $\delta$ such that $C_{w,n,\alpha}'\delta^{1/q}<1/8$. Then by \eqref{estimate of low part of fractional}, for any $x\in Q_j$, there exists $j_0\geq1$ such that for any $j>j_0$,
\begin{align*}
jI_\alpha(f_j^{low})(x)\leq C_{w,n,\alpha}'(\delta^{1/q}+je^{-\gamma_w4^j/(2q)})<1/4.
\end{align*}
This proves \eqref{case 2 claim}.

Hence, by \eqref{case 2 claim} and the weighted endpoint estimate of $I_\alpha$, for any $j>j_0$,
\begin{align}\label{transform to high}
&w(\{x\in\mathbb{R}^n:j\chi_{Q_j}(x)I_\alpha(|f|)(x)>1/2\})\\
&\quad\leq w(\{x\in\mathbb{R}^n:\chi_{Q_j}(x)I_\alpha(f_j^{high})(x)>1/(4j)\})\nonumber\\
&\quad\leq C_w\Big(4j\int_{\mathbb{R}^n}f_j^{high}(x)w(x)^{1/q}dx\Big)^q\nonumber\\
&\quad=C_w\Big(4j\int_{\mathbb{R}^n}\frac{\Phi(f_j^{high}(x))}{\log(e+f_j^{high}(x))}w(x)^{1/q}dx\Big)^q
\nonumber.
\end{align}
Note that there exists $j_1\geq1$ such that for any $j>j_1$,
\begin{align}\label{log estimate2}
\log(e+f_j^{high})&\geq\log(f_j^{high})\\
&\geq\log\Big(\delta j^{-q}\Big(\int_{\mathbb{R}^n}\Phi(|f|)w^{1/q}\Big)^{1-q}\Big)\nonumber\\
&\geq\log(\delta j^{-q}e^{-4^j\gamma_w(1-q)/(2q)})\nonumber\\
&=\log\delta-q\log j+4^j\gamma_w(q-1)/(2q)\geq4^j\gamma_w(q-1)/(4q).\nonumber
\end{align}
Inserting \eqref{log estimate2} into \eqref{transform to high}, for any $j>\max\{j_1,j_0\}$, we further have
\begin{align*}
w(\{x\in\mathbb{R}^n:j\chi_{Q_j}(x)I_\alpha(|f|)(x)>1/2\})\leq C_{w,n,\alpha}\Big(j4^{-j}\int_{\mathbb{R}^n}\Phi(|f(x)|)w(x)^{1/q}dx\Big)^q.
\end{align*}

Finally, to prove the desired result for all $j\geq1$, we need to show the conclusion holds for $1\leq j\leq\max\{j_1,j_0\}$. As a consequence of weighted endpoint estimate of $I_\alpha$,
\begin{align*}
&w(\{x\in\mathbb{R}^n:j\chi_{Q_j}(x)I_\alpha(|f|)(x)>1/2\})\\
&\quad\leq C_w\Big(2j\int_{\mathbb{R}^n}|f(x)|w(x)^{1/q}dx\Big)^q\\
&\quad\leq C_w\Big(2j4^j4^{-j}\int_{\mathbb{R}^n}|\Phi(f(x))|w(x)^{1/q}dx\Big)^q\\
&\quad\leq C_w4^{\max\{j_1,j_0\}}\Big(j4^{-j}\int_{\mathbb{R}^n}\Phi(|f(x)|)w(x)^{1/q}dx\Big)^q.
\end{align*}
This completes the proof of Lemma \ref{key lemma2}.
\end{proof}

\begin{proof}[Proof of Theorem \ref{counterexample}]
Let $Q_0:=[-1/2,1/2]^n$. Let $\{Q_j\}_j$ be a collection of disjoint cubes with $l(Q_j)=e^{-4^j}$, that is contained in $Q_0$. Define $b(x)=\sum_{j=1}^{\infty}j\chi_{Q_j}(x)$. It is not hard to verify that $b\not\in {\rm BMO}(\mathbb{R}^{n})$. In fact, let $\{\widetilde{Q_j}\}_j$ be a family of disjoint cubes that satisfy $\widetilde{Q_j}=2^{1/n}Q_j$. Then $b_{\widetilde{Q_j}}=j/2$ and
\begin{align*}
&\frac{1}{|\widetilde{Q_j}|}\int_{\widetilde{Q_j}}|b(x)-b_{\widetilde{Q_j}}|dx\\
&\quad=\frac{1}{|\widetilde{Q_j}|}\int_{\widetilde{Q_j}\backslash Q_j}|b(x)-b_{\widetilde{Q_j}}|dx+\frac{1}{|\widetilde{Q_j}|}\int_{Q_j}|b(x)-b_{\widetilde{Q_j}}|dx\\
&\quad=\frac{j}{2}\frac{|\widetilde{Q_j}\backslash Q_j|+|Q_j|}{|\widetilde{Q_j}|}=\frac{j}{2}.
\end{align*}
This proves $b\not\in {\rm BMO}(\mathbb{R}^{n})$.

For such a function $b$, we shall prove that for any $t>0$ and $w\in A_1$,
\begin{align*}
w(\{x\in\mathbb{R}^n:[b,I_\alpha]f(x)>t\})\leq C_{w,n,\alpha}\Big(\int_{\mathbb{R}^n}\Phi\Big(\frac{|f(x)|}{t}\Big)w(x)^{1/q}dx\Big)^q.
\end{align*}
If we have done, the conclusion holds as $\Theta(w)\geq w^{1/q}$ and $\Psi(t)\geq t^q$.

By homogeneous, it suffices to prove
\begin{align*}
w(\{x\in\mathbb{R}^n:[b,I_\alpha]f(x)>1\})\leq C_{w,n,\alpha}\Big(\int_{\mathbb{R}^n}\Phi(|f(x)|)w(x)^{1/q}dx\Big)^q.
\end{align*}
Since $\sum_{j=1}^\infty2^{-j}=1$,
\begin{align}\label{transform to sum}
w(\{x\in\mathbb{R}^n:[b,I_\alpha]f(x)>1\})\leq \sum_{j=1}^\infty w(\{x\in\mathbb{R}^n:[j\chi_{Q_j},I_\alpha]f(x)>2^{-j}\}).
\end{align}
Thanks to Lemmas \ref{key lemma1} and \ref{key lemma2}, we get
\begin{align*}
w(\{x\in\mathbb{R}^n:[j\chi_{Q_j},I_\alpha]f(x)>1\})\leq C_{w,n,\alpha}\Big(j4^{-j}\int_{\mathbb{R}^n}\Phi(|f(x)|)w(x)^{1/q}dx\Big)^q.
\end{align*}
From this and \eqref{transform to sum}, we obtain
\begin{align*}
&w(\{x\in\mathbb{R}^n:[j\chi_{Q_j},I_\alpha]f(x)>1\})\\
&\quad\leq C_{w,n,\alpha}\sum_{j=1}^\infty\Big(j4^{-j}\int_{\mathbb{R}^n}\Phi(2^j|f(x)|)w(x)^{1/q}dx\Big)^q\\
&\quad\leq C_{w,n,\alpha}\sum_{j=1}^\infty\Big(j4^{-j}(j+1)2^j\int_{\mathbb{R}^n}\Phi(|f(x)|)w(x)^{1/q}dx\Big)^q\\
&\quad\leq C_{w,n,\alpha}\Big(\int_{\mathbb{R}^n}\Phi(|f(x)|)w(x)^{1/q}dx\Big)^q,
\end{align*}
as we wanted.
\end{proof}

\section{A counterexample to two-weighted endpoint estimate of Bloom type for commutator of Hilbert transform}
In this section, we give a counterexample to illustrate that two-weighted endpoint estimate of Bloom type for commutator of Hilbert transform does not hold.

\begin{lemma}\label{weight belong to A1}
Let $\lambda(x)=(1+|x|)^{-1/2}$, $x\in\mathbb{R}$. Then $\lambda\in A_1$.
\end{lemma}
\begin{proof}
Fix an interval $I:=(a,b)\subset\mathbb{R}$, we need to prove
\begin{align*}
\frac{1}{|I|}\int_I\lambda(x)dx\lesssim\essi_{x\in I}\lambda(x).
\end{align*}
We consider two cases: $\inf_{x\in I}\lambda(x)\geq b-a$ and $\inf_{x\in I}\lambda(x)<b-a$.

\textbf{Case 1. $\inf_{x\in I}|x|\geq b-a$}:\quad By a simple observation, for any $x\in I$,
\begin{align}\label{sup controled by inf}
|x|\leq\inf_{x\in I}|x|+(b-a)\leq2\inf_{x\in I}|x|.
\end{align}
It follows that
\begin{align*}
\essi_{x\in I}\lambda(x)=\Big(1+\sup_{x\in I}|x|\Big)^{-1/2}\geq\Big(1+2\inf_{x\in I}|x|\Big)^{-1/2}\geq2^{-1/2}\inf_{x\in I}|x|.
\end{align*}
So we have
\begin{align*}
\frac{1}{|I|}\int_I\lambda(x)dx&=\frac{1}{|I|}\int_I(1+|x|)^{-1/2}dx\\
&\leq\Big(1+\inf_{x\in I}|x|\Big)^{-1/2}\leq\Big(\inf_{x\in I}|x|\Big)^{-1/2}\\
&\leq\sqrt{2}\essi_{x\in I}\lambda(x).
\end{align*}

\textbf{Case 2. $\inf_{x\in I}|x|<b-a$}:\quad Again by \eqref{sup controled by inf}, for any $x\in I$, $|x|\leq2(b-a)$. Then
\begin{align*}
\frac{1}{|I|}\int_I\lambda(x)dx&\leq\frac{1}{b-a}\int_{-\sup_{x\in I }|x|}^{\sup_{x\in I}|x|}(1+|x|)^{-1/2}dx\\
&=\frac{2}{b-a}\int_{0}^{\sup_{x\in I}|x|}(1+|x|)^{-1/2}dx\\
&=\frac{4}{b-a}\Big[\Big(1+\sup_{x\in I}|x|\Big)^{1/2}-1\Big]\\
&=\frac{4}{b-a}\frac{\sup_{x\in I}|x|}{1+\Big(1+\sup_{x\in I}|x|\Big)^{1/2}}\\
&\leq8\Big(1+\sup_{x\in I}|x|\Big)^{-1/2}=8\essi_{x\in I}\lambda(x).
\end{align*}

Combined with Case 1 and Case 2, we conclude that $\lambda\in A_1$.
\end{proof}

\begin{lemma}\label{g belong to BMO}
Let $y\in\mathbb{R}$ and
\[
g(y)\\
:=
\begin{cases}
\log\frac{1}{|y|},& 0<|y|<1,\\
0,& y=0,\\
0,& |y|\geq1.
\end{cases}
\]
Then $g\in {\rm BMO(\mathbb{R})}$.
\end{lemma}
\begin{proof}
Note that we can rewrite
\[
g(y)\\
:=
\begin{cases}
\max\{-\log|y|,0\}, & y\neq0,\\
0,& y=0.
\end{cases}
\]
By the fact that $\log|y|\in {\rm BMO(\mathbb{R})}$, see \cite[Example 3.13]{Gra1}, and
\begin{equation*}
\max\{-\log|y|,0\}\leq3\|-\log|u|\|_{\rm BMO(\mathbb{R})}/2,
\end{equation*}
see \cite[Proposition 3.1.2]{Gra1}, we get the desired result.
\end{proof}

\begin{lemma}\label{b belongs to weighted BMO}
Let $\eta(x)=(1+|x|)^{1/2}$. Let $x\in\mathbb{R}$ and $b(x)=\sum_{j=1}^\infty2^jg(x-4^j)$, where
\[
g(y)\\
:=
\begin{cases}
\log\frac{1}{|y|},& 0<|y|<1,\\
0,& y=0,\\
0,& |y|\geq1.
\end{cases}
\]
Then $b\in {\rm BMO_\eta(\mathbb{R})}$.
\end{lemma}
\begin{proof}
Fix an interval $I:=(a,b)\subset\mathbb{R}$. If $I\cap(4^j-1,4^j+1)=\emptyset$ for any $j\in\mathbb{N}$, then $b=0$, it is trivial that $b\in {\rm BMO_\eta(\mathbb{R})}$. Hence, it suffices toc consider $I\cap(4^j-1,4^j+1)\neq\emptyset$. Let $J_i:=(4^i-2,4^i+2)$, $i\in\mathbb{N}$. Then $J_i\cap J_l=\emptyset$, $i\neq l$, and for $x\in J_i$,
\begin{align}\label{bizhi weight estimate}
\eta(x)\geq(4^i-1)^{1/2}\geq(4^i/2)^{1/2}.
\end{align}
We split the proof of this lemma into two cases: $b-a\leq1$ and $b-a>1$.

\textbf{Case 1. $b-a\leq1$}:\quad  There exists only one $j_0\in\mathbb{N}$ such that
\begin{align*}
I\cap(4^{j_0}-1,4^{j_0}+1)\neq\emptyset
\end{align*}
as $b-a\leq1$. Furthermore, $I\subset J_{j_0}$. Based on these facts, we have $b(x)=2^{j_0}g(x-4^{j_0})$ and
\begin{align*}
\int_{I}|b(x)-b_I|dx&=2^{j_0}\int_{I}|g(x-4^{j_0})-(g(\cdot-4^{j_0}))_I|dx\\
&\leq C_g2^{j_0}|I|\leq \sqrt{2}C_g\eta(I),
\end{align*}
where we use Lemma \ref{g belong to BMO} and \eqref{bizhi weight estimate} in the penultimate inequality and the last inequality, respectively.

\textbf{Case 2. $b-a>1$}:\quad Denote
\begin{align*}
L(I):=\{j\in\mathbb{N}:I\cap(4^j-1,4^j+1)\neq\emptyset\}.
\end{align*}
For each $j\in L(I)$, we claim that $|I\cap J_j|\geq1$. In fact, take $y_j\in I\cap(4^j-1,4^j+1)$, then from $4^j-1<y_j<4^j+1$, we have
\begin{align}\label{yj}
y_j-(4^j-2)>1,~(4^j+2)-y_j>1.
\end{align}
If $I\subset J_j$,
$$|I\cap J_j|=|I|=b-a>1.$$
If $I\not\subset J_j$, \eqref{yj} indicates that $I$ has at least one endpoint outside $J_j$. Assume $a\leq 4^j-2$, then
\begin{align*}
[4^j-2,y_j]\subset I\cap J_j,
\end{align*}
it follows that $|I\cap J_j|\geq y_j-(4^j-2)>1$ by \eqref{yj}. The same result can be obtained in case $b\geq 4^{j}+2$. This proves the claim.

By the claim $|I\cap J_j|\geq1$ and \eqref{bizhi weight estimate}, we deduce that
\begin{align*}
\int_I\eta(x)dx&\geq\sum_{j\in L(I)}\int_{I\cap J_j}\eta(x)dx\\
&\geq\frac{\sqrt{2}}{2}\sum_{j\in L(I)}2^j|I\cap J_j|\\
&\geq\frac{\sqrt{2}}{2}\sum_{j\in L(I)}2^j.
\end{align*}
From this, we conclude that
\begin{align*}
\int_I|b(x)-b_I|dx&\leq2\int_I|b(x)|dx\\
&=2\int_Ib(x)dx\\
&\leq2\sum_{j\in L(I)}2^j\int_{-1}^{1}g(x)dx\\
&=4\int_0^1\log\frac{1}{|x|}dx\sum_{j\in L(I)}2^j\\
&=4\sum_{j\in L(I)}2^j\leq4\sqrt{2}\eta(I).
\end{align*}
Therefore, we obtain $b\in {\rm BMO_\eta(\mathbb{R})}$.
\end{proof}

Now, we are in the position to prove Theorem \ref{counterexample1}.
\begin{proof}[Proof of Theorem \ref{counterexample1}]
Assume that for $b\in{\rm BMO}_\eta(\mathbb{R})$,
\begin{align}\label{contraction}
\lambda(\{x\in\mathbb{R}:|[b,H]f(x)|>t\})\leq C\int_{\mathbb{R}}\Phi\Big(\frac{|f(x)|}{t}\Big)\mu(x)dx
\end{align}
holds for any $\mu,\lambda\in A_1$. Let $\mu=1$, $\lambda$ and $b$ be as in Lemma \ref{weight belong to A1} and Lemma \ref{b belongs to weighted BMO}, respectively. Let $E_j:=(4^j+2^{-j-1},4^j+2^{-j})$, $F_j:=(4^j+1,4^j+2)$. Then for $x\in F_j$, $y\in E_j$,
\begin{align*}
b(x)=0,~b(y)\geq2^j\log2^j,~0<x-y<2,
\end{align*}
which implies that for any $x\in F_j$,
\begin{align*}
|[b,H](\chi_{E_j})(x)|&=\Big|\int_{E_j}\frac{b(x)-b(y)}{x-y}dy\Big|\\
&=\int_{E_j}\frac{b(y)}{x-y}dy\\
&\geq2^{-1}2^{-j-1}2^j\log2^j=j\log2/4.
\end{align*}
Hence,
$$F_j\subset\{x\in\mathbb{R}:|[b,H](\chi_{E_j})(x)|>j\log2/8\}.$$
By \eqref{contraction}, we have
\begin{align*}
(1+4^j+2)^{-1/2}|F_j|&\leq\lambda(F_j)\\
&\leq\lambda(\{x\in\mathbb{R}:|[b,H](\chi_{E_j})(x)|>j\log2/8\})\\
&\leq c\int_{E_j}\Phi\Big(\frac{8}{j\log2}\Big)dx\\
&=c 2^{-j-1}\Phi\Big(\frac{8}{j\log2}\Big).
\end{align*}
This obviously implies that $1\leq c\Phi\Big(\frac{8}{j\log2}\Big)$. It is a contraction when $j\rightarrow\infty$. So we finish the proof of Theorem \ref{counterexample1}.
\end{proof}
\subsection*{Declaration of competing interest}
The authors declare that there is no competing interest.

\subsection*{Acknowledgements}
The authors acknowledge the use of AI tools during the exploratory stage of two-weighted endpoint estimate of Bloom type for commutator of Hilbert transform. All other mathematical contents are solely due to the authors.
%%%%%%%%%%%% References %%%%%%%%%%%%%
%%
%<Author name> is written as Initial of Given Name, and Family Name.
%<Title> is written in roman letters.
%<Journal name> should be abbreviated according to
% the MR Serials Abbreviations List of Mathematical Reviews:
% (Abbreviations of Names of Serials; http://www.ams.org/mr-database)
%For <Pages>, use en-dash "--" between page numbers.
%%


\begin{thebibliography}{99}

%%% Journals %%%

%<Author name>, <Title>, <Journal name> <Volume> (<Year>), <Pages>.
\bibitem{A}N. Accomazzo,
A characterization of BMO in terms of endpoint bounds for commutators of singular integrals,
Israel J. Math. 228(2) (2018), 787--800.

\bibitem{AMR}N. Accomazzo, J.C. Mart\'{i}nez-Perales I.P. Rivera-R\'{i}os,
On Bloom-type estimates for iterated commutators of fractional integrals,
Indiana Univ. Math. J. 69(4) (2020), 1207--1230.

\bibitem{AJS}M. Akcoglu, R.L. Jones and P. Schwartz,
Variation in probability, ergodic theory and analysis,
Illinois J. Math. 42 (1998), 154--177.

\bibitem{BFHR}J.J. Betancor, J.C. Farina, E. Harbour and L. Rodriguez-Mesa,
$L^p$-boundedness properities of variation operators in the Schr\"odinger setting,
Rev. Mat. Complut. 26(2) (2013), 485--534.

\bibitem{Bloom}S. Bloom,
A commutator theorem and weighted BMO,
Trans. Amer. Math. Soc. 292(1) (1985), 103--122.

\bibitem{Bo}J. Bourgain,
Pointwise ergodic theorems for arithmetric sets,
Publ. Math. Inst. Hautes \'{E}tudes Sci. 69(1) (1989), 5--45.

\bibitem{CaJRW1}J.T. Campbell, R.L. Jones, K. Reinhold and M. Wierdl,
Oscillations and variation for the Hilbert transform,
Duke Math. J. 105(1) (2000), 59--83.

\bibitem{CaJRW2}J.T. Campbell, R.L. Jones, K. Reinhold and M. Wierdl,
Oscillations and variation for singular integrals in higher dimensions,
Trans. Amer. Math. Soc. 355(5) (2003), 2115--2137.

\bibitem{ChWiWo}S.-Y. A. Chang, J.M. Wilson and T.H. Wolff,
Some weighted norm inequalities concerning the Schr\"{o}dinger operators,
Comment. Math. Helv. 60(2) (1985), 217--246.

\bibitem{Ch}S. Chanillo,
A note on commutators,
Indiana Univ. Math. J. 31 (1982), 7--16.

\bibitem{ChDiHoLi}Y. Chen, Y. Ding, G. Hong and H. Liu,
Variational inequalities for the commutators of rough operators with BMO functions,
Sci. China Math. 64(11) (2021), 2437--2460.

\bibitem{CoRoWe}R.R. Coifman, R. Rochberg and G. Weiss,
Factorization theorems for Hardy spaces in several variables,
Ann. of Math. (2) 103(3) (1976), 611--635.

\bibitem{CrFi2003}D. Cruz-Uribe and A. Fiorenza,
Endpoint estimates and weighted norm inequalities for commutators of fractional integrals,
Publ. Mat. 47(1) (2003), 103--131.

\bibitem{CrFi2007}D. Cruz-Uribe and A. Fiorenza,
Weighted endpoint estimates for commutators of fractional integrals,
Czechoslovak Math. J. 57(132)(1) (2007), 153--160.

\bibitem{DiLu1999}Y. Ding and S. Lu,
Higher order commutators for a class of rough operators,
Ark. Mat. 37(1) (1999), 33--44.

\bibitem{DiLuZh}Y. Ding, S. Lu and P. Zhang,
Weak estimates for commutators of fractional integral operators,
Sci. China Ser. A 44(7) (2001), 877--888.

\bibitem{DiXu0}Y. Ding and Q. Xue,
Weak type (1,1) bounds for a class of the Littlewood-Paley operators,
J. Math. Soc. Japan 57(1) (2005), 183--194.

\bibitem{DiXu}Y. Ding and Q. Xue,
Endpoint estimates for commutators of a class of Littlewood-Paley operators,
Hokkaido Math. J. 36(2) (2007), 245--282.

\bibitem{DuGoKuLiWiYa}X.T. Duong, R. Gong, M.-J. S. Kuffner, J. Li, B.D. Wick and D. Yang,
Two weight commutators on spaces of homogeneous type and applications,
J. Geom. Anal. 31(1) (2021), 980--1038.

\bibitem{Gra}L. Grafakos,
Classical Fourier analysis,
Third edition. Graduate Texts in Mathematics,
Springer, New York, (2014), pp. xviii+638. ISBN: 978-1-4939-1194-3.

\bibitem{Gra1}L. Grafakos,
Modern Fourier analysis,
Third edition. Graduate Texts in Mathematics, 250.
Springer, New York, 2014. xvi+624 pp. ISBN: 978-1-4939-1229-2; 978-1-4939-1230-8.


\bibitem{GLW}W. Guo, J. Lian and H. Wu,
The unified theory for the necessity of bounded commutators and applications,
J. Geom. Anal. 30(4) (2020), 3995--4035.

\bibitem{GuoWenWuYang}W. Guo, Y. Wen, H. Wu and D. Yang,
Variational characterizations of weighted Hardy spaces and weighted BMO spaces,
Proc. Roy. Soc. Edinburgh Sect. A 152(6) (2022), 1613--1632.

\bibitem{HLW}I. Holmes, M. T. Lacey and B. D. Wick,
Commutators in the two-weight setting,
Math. Ann. 367(1-2) (2016), 51--80.

\bibitem{HuLiWa}B. Hu, J. Li and N.A. Wagner,
Boundedness and compactness of Bergman projection commutators in two-weight setting,
J. Fourier Anal. Appl. 32(1) (2026), Paper No. 15, 39pp.

\bibitem{HuTa}G. Hu and X. Tao,
An endpoint estimate for the commutators of singular integral operators with rough kernels,
Potential Anal. 58(2) (2023), 241--262.

\bibitem{Hy}T.P. Hyt\"{o}nen,
The $L^p$-to-$L^q$ boundedness of commutators with applications to the Jacobian operator,
J. Math. Pures Appl. (9) 156 (2021), 351--391.

\bibitem{Ja}S. Janson,
Mean oscillation and commutators of singular integral operators,
Ark. Mat. 16(2) (1978), 263--270.

\bibitem{Jo}R.L. Jones,
Ergodic theory and connections with analysis and probability,
New York J. Math. 3A (1997/98), Proceedings of the New York Journal of Mathematics Conference, 31--67.

\bibitem{JR}R.L. Jones and K. Reinhold,
Oscillation and variation inequalities for convolution powers,
Ergodic Theory Dynam. Sys. 21(6) (2001), 1809--1829.

\bibitem{Ke}C. Kenig,
Harmonic analysis techniques for second order elliptic boundary value problems,
CBMS Regional Conference Series in Mathematics, 83. Published for the Conference Board of the Mathematical Sciences, Washington, DC; by the American Mathematical Society, Providence, RI, 1994. xii+146pp. ISBN: 0-8218-0309-3.

\bibitem{LaHuTa}J. Lan, X. Tao and G. Hu,
Weak type endpoint estimates for the commutators of rough singular integral operators,
Math. Inequal. Appl. 23(4) (2020), 1179--1195.

\bibitem{Le}D. L\'{e}pingle,
La variation d'order $p$ des semi-martingales,
Z. Wahrscheinlichkeitstheorie und Verw. Gebiete, 36(4) (1976), 295--316.

\bibitem{LeOmRi2017}A.K. Lerner, S. Ombrosi and I.P. Rivera-R\'{i}os,
On pointwise and weighted estimates for commutators of Calder\'{o}n-Zygmund operators,
Adv. Math. 319 (2017), 153--181.

\bibitem{LeOmRi}A.K. Lerner, S. Ombrosi and I.P. Rivera-R\'{i}os,
Commutators of singular integrals revisited,
Bull. Lond. Math. Soc. 51(1) (2019), 107--119.

\bibitem{LiLi}J. Li and K. Li,
Complex meadian method,
arXiv: 2609.28088v2.

\bibitem{Li}K. Li,
Multilinear commutators in the two-weight setting,
Bull. Lond. Math. Soc. 54(2) (2022), 568--589.

\bibitem{LiMaVu}K. Li, H. Martikainen and E. Vuorinen,
Bloom-type inequality for bi-parameter singular integrals: efficient proof and iterated commutators,
Int. Math. Res. Not. IMRN 2021(11) (2021), 8153--8187.


\bibitem{LW}F. Liu and H. Wu,
A criterion on oscillation and variation for the commutators of singular integrals,
Forum Math. 27(1) (2015), 77--97.

\bibitem{MuWh}B. Muckenhoupt and R. Wheeden,
Weighted norm inequalities for fractional integrals,
Trans. Amer. Math. Soc. 192 (1974), 261--274.

\bibitem{P}C. P\'{e}rez,
Endpoint estimates for commutators of singular integral operators,
J. Funct. Anal. 128(1) (1995), 163--185.

\bibitem{SaYa}M. Sakamoto and K. Yabuta,
Boundedness of Marcinkiewicz functions,
Studia Math. 135(2) (1999), 103--142.

\bibitem{ST}C. Segovia and J.L. Torrea,
Weighted inequalities for commutators of fractional and singular integrals, in: Conference on Mathematical Analysis (El Escorial, 1989),
Publ. Mat. 35 (1991), 209--235.

\bibitem{ST1993}C. Segovia and J.L. Torrea,
Higher order commutators for vector-valued Calder¨®n¨CZygmund operators,
Trans. Amer. Math. Soc. 336(2) (1993), 537--556.

\bibitem{St}E.M. Stein,
On the function of Littlewood-Paley, Lusin and Marcinkiewicz,
Trans. Amer. Math. Soc., 88 (1958), 430--466.

\bibitem{TW}A. Torchinsky and S. Wang,
A note on the Marcinkiewicz integral,
Colloq. Math. 60/61(1) (1990), 235--243.

\bibitem{Uc}A. Uchiyama,
On the compactness of operators of Hankel type,
Tohoku Math. J. (2) 30(1) (1978),163--171.

\bibitem{WenWuXue}Y. Wen, H. Wu and Q. Xue
Sparse domination and weighted inequalities for the $\rho$-variation of singular integrals and commutators,
J. Geom. Anal. 32(12) (2022), Paper No. 297, 30pp.

\bibitem{Zhang}P. Zhang,
Weighted endpoint estimates for commutators of Marcinkiewicz integrals,
Acta Math. Sin. (Engl. Ser.) 24(8) (2008), 1387--1400.
\end{thebibliography}
\end{document}